\documentclass{gtart}
\def\ifplaintex{\expandafter\ifx\csname documentclass\endcsname\relax}

\ifplaintex 
\else
\fi

\expandafter\ifx\csname beginpicture\endcsname\relax
\expandafter\ifx\csname documentclass\endcsname\relax
\input pictex \else
\input prepictex \input pictex \input postpictex \fi\fi

\def\gt{{\mathsurround=0pt\it $\cal G\mskip-2mu$eometry \&\ 
$\cal T\!\!$opology}}        %  journal title in recommended style

\def\gtp{{\mathsurround=0pt\it $\cal G\mskip-2mu$eometry \&\ 
$\cal T\!\!$opology $\cal P\!$ublications}}  % GT publications

\def\lognumber#1{\def\thelognumber{#1}}
\def\volumenumber#1{\def\thevolumenumber{#1}}
\def\papernumber#1{\def\thepapernumber{#1}}
\def\volumeyear#1{\def\thevolumeyear{#1}}

\def\pagenumbers#1#2{\def\startpage{#1}\def\finishpage{#2}}
\def\published#1{\def\publishdate{#1}}
\def\proposed#1{\def\theproposer{#1}}
\def\seconded#1{\def\theseconders{#1}}
\def\received#1{\def\receiveddate{#1}}
\def\revised#1{\def\reviseddate{#1}}
\def\accepted#1{\def\accepteddate{#1}}

\long\def\asciiabstract#1{\long\def\theasciiabstract{#1}}

\let\\\par\let\thelognumber\relax
\let\thevolumenumber\relax\let\thepapernumber\relax
\let\thevolumeyear\relax\let\thesamplenumber\relax\let\startpage\relax
\let\finishpage\relax\let\publishdate\relax\let\receiveddate\relax
\let\reviseddate\relax\let\accepteddate\relax\let\theasciititle\relax
\let\theasciiauthors\relax
\let\theasciiabstract\relax
\let\theasciiemail\relax\let\theshortauthors\relax\let\theshorttitle\relax

\long\def\maketitlep{   % start of definition of \maketitlep

\count0=\startpage

\gt\hfill      %   Journal title (top left) 
\beginpicture
\setcoordinatesystem units <0.33truein, 0.33truein> point at 2.2 0.9
\setplotsymbol ({$\cal G$})
\plotsymbolspacing=9truept
\circulararc 315 degrees from 0 1 center at 0 0
\setplotsymbol ({$\cal T$})
\circulararc 315 degrees from 1 -1 center at 1 0
\endpicture
\break
{\small\ifx\thesamplenumber\relax % sample?  
Volume \else Sample
\fi\thevolumenumber\ (\thevolumeyear)
\startpage--\finishpage\nl
Published: \publishdate}
\vglue 0.5truein plus 0.4fil minus 0.1truein

{\parskip=0pt\leftskip 0pt plus 1fil\def\\{\par\smallskip}{\ifplaintex\large
\else\Large\fi\bf\thetitle}\par\medskip}   

\vglue 0pt plus 0.1fil 

{\parskip=0pt\leftskip 0pt plus 1fil\def\\{\par}{\sc\theauthors}
\par\medskip}

\vglue 0pt plus 0.1fil 

{\small\parskip=0pt\let\newline\\
{\leftskip 0pt plus 1fil\def\\{\par}{\sl\theaddress}\par}
\expandafter\ifx\theemail\relax    % email address?
\relax\else\vglue 5pt plus 0.02fil minus 2pt\def\\{\stdspace{\rm 
and}\stdspace} 
\cl{Email:\stdspace\tt\theemail}\fi
\ifx\theurl\relax                  % URL given?
\relax\else\vglue 5pt plus 0.02fil minus 2pt\def\\{\stdspace{\rm 
and}\stdspace}
\cl{URL:\stdspace\tt\theurl}\fi\par}

\vglue 7pt plus 0.3fil minus 3pt

{\bf Abstract}
\vglue 5pt plus 0.1fil minus 2pt

\theabstract

\vglue 7pt plus 0.3fil minus 3pt

{\bf AMS Classification numbers}\quad Primary:\quad \theprimaryclass

Secondary:\quad \thesecondaryclass

\vglue 5pt plus 0.3fil minus 2pt

{\bf Keywords}\quad \thekeywords

\vglue 10pt plus 0.5fil minus 5pt

{\small  Proposed: \theproposer\hfill Received: \receiveddate\nl
Seconded: \theseconders\hfill 
\ifx\reviseddate\relax                         % paper revised?
Accepted: \accepteddate                        % no
\else
Revised: \reviseddate                          % yes
\fi}
\eject
}       %  end of definition of \maketitlep

\let\maketitlepage\maketitlep
\let\maketitle\maketitlepage

\font\phead=cmsl9 scaled 950
\font=cmsl9 scaled 1050
\font\pnum=cmbx10 scaled 913
\font\lnum=cmbx10 
\font\pfoot=cmsl9 scaled 950
\font=cmsl9 scaled 1050
\ifplaintex
\headline{\vbox to 0pt{\vskip -4.5mm\line{\small\phead\ifnum
\count0=\startpage ISSN 1364-0380 (on line)
1465-3060 (printed) \hfill {\pnum\folio}\else\ifodd\count0\def\\{ }% 
\ifx\theshorttitle\relax\thetitle\else\theshorttitle\fi\hfill{\pnum\folio}
\else\def\\{ and }{\pnum\folio}\hfill\ifx\theshortauthors\relax\theauthors
\else\theshortauthors\fi\fi\fi}\vss}}
\footline{\vbox to 0pt{\vglue 0mm\line{\small\pfoot\ifnum\count0=\startpage
\copyright\ \gtp\hfill\else
\gt, Volume \thevolumenumber\ (\thevolumeyear)\hfill\fi}\vss
}}
\else
\makeatletter
\def\@oddhead{{\small\lhead\ifnum\count0=\startpage ISSN 1364-0380 (on line)
1465-3060 (printed) \hfill {\lnum\number\count0}\else\ifodd\count0
\def\\{ }\ifx\theshorttitle\relax \thetitle \else\theshorttitle\fi\hfill
{\lnum\number\count0}\else\def\\{ and }{\lnum\number\count0}
\hfill\ifx\theshortauthors\relax 
\theauthors\else\theshortauthors\fi\fi\fi}}\def\@evenhead{\@oddhead}
\def\@oddfoot{\small\lfoot\ifnum\count0=\startpage\copyright\ \gtp\hfill\else
\gt, Volume \thevolumenumber\ (\thevolumeyear)\hfill\fi}
\def\@evenfoot{\@oddfoot}
\makeatother
\fi

\newwrite\gtoutfile
\long\gdef\makeheadfile{  %%% start of definition of \makeheadfile
{\def\\{, }\def\s{ }
\immediate\openout\gtoutfile head.xxx
\immediate\write\gtoutfile{To: math@arxiv.org}
\immediate\write\gtoutfile{Subject: put or rep NNNNN:pppp}
\immediate\write\gtoutfile{--text follows this line--}
\immediate\write\gtoutfile{Proxy-for: \ifx\theasciiauthors\relax
\theauthors\else\theasciiauthors\fi\s<\ifx\theasciiemail\relax\theemail\else\theasciiemail\fi>}
\immediate\write\gtoutfile{\noexpand\\}
\immediate\write\gtoutfile{Authors: \ifx\theasciiauthors\relax
\theauthors\else\theasciiauthors\fi}
\immediate\write\gtoutfile{Title: \ifx\theasciititle\relax
\thetitle\else\theasciititle\fi}
\immediate\write\gtoutfile{Subj-class: GT or SG or MG etc}
\immediate\write\gtoutfile{MSC-class: \theprimaryclass\ifx\thesecondaryclass\relax\else, \thesecondaryclass\fi}
\immediate\write\gtoutfile{Journal-ref: Geom. Topol. \thevolumenumber
(\thevolumeyear) \startpage-\finishpage}
\immediate\write\gtoutfile{Comments: Published by Geometry and Topology at}
\immediate\write\gtoutfile{\s\s http://www.maths.warwick.ac.uk/gt/GTVol\thevolumenumber/paper\thepapernumber.abs.html}
\immediate\write\gtoutfile{\noexpand\\}
\immediate\write\gtoutfile{}
\ifx\theasciiabstract\relax
\immediate\write\gtoutfile{\theabstract}\else
\immediate\write\gtoutfile{\theasciiabstract}\fi
\immediate\write\gtoutfile{}
\immediate\write\gtoutfile{\noexpand\\}
\immediate\write\gtoutfile{}
\immediate\closeout\gtoutfile}}  %%% end of definition of \makeheadfile

\def\maketitlepage{\maketitlep\makeheadfile}
\let\maketitle\maketitlepage

\def\ifplaintex{\expandafter\ifx\csname documentclass\endcsname\relax}

\ifplaintex 
\else
\fi

\expandafter\ifx\csname beginpicture\endcsname\relax
\expandafter\ifx\csname documentclass\endcsname\relax
\input pictex \else
\input prepictex \input pictex \input postpictex \fi\fi

\def\gt{{\mathsurround=0pt\it $\cal G\mskip-2mu$eometry \&\ 
$\cal T\!\!$opology}}        %  journal title in recommended style

\def\gtp{{\mathsurround=0pt\it $\cal G\mskip-2mu$eometry \&\ 
$\cal T\!\!$opology $\cal P\!$ublications}}  % GT publications

\def\lognumber#1{\def\thelognumber{#1}}
\def\volumenumber#1{\def\thevolumenumber{#1}}
\def\papernumber#1{\def\thepapernumber{#1}}
\def\volumeyear#1{\def\thevolumeyear{#1}}

\def\pagenumbers#1#2{\def\startpage{#1}\def\finishpage{#2}}
\def\published#1{\def\publishdate{#1}}
\def\proposed#1{\def\theproposer{#1}}
\def\seconded#1{\def\theseconders{#1}}
\def\received#1{\def\receiveddate{#1}}
\def\revised#1{\def\reviseddate{#1}}
\def\accepted#1{\def\accepteddate{#1}}

\long\def\asciiabstract#1{\long\def\theasciiabstract{#1}}

\let\\\par\let\thelognumber\relax
\let\thevolumenumber\relax\let\thepapernumber\relax
\let\thevolumeyear\relax\let\thesamplenumber\relax\let\startpage\relax
\let\finishpage\relax\let\publishdate\relax\let\receiveddate\relax
\let\reviseddate\relax\let\accepteddate\relax\let\theasciititle\relax
\let\theasciiauthors\relax
\let\theasciiabstract\relax
\let\theasciiemail\relax\let\theshortauthors\relax\let\theshorttitle\relax

\long\def\maketitlep{   % start of definition of \maketitlep

\count0=\startpage

\gt\hfill      %   Journal title (top left) 
\beginpicture
\setcoordinatesystem units <0.33truein, 0.33truein> point at 2.2 0.9
\setplotsymbol ({$\cal G$})
\plotsymbolspacing=9truept
\circulararc 315 degrees from 0 1 center at 0 0
\setplotsymbol ({$\cal T$})
\circulararc 315 degrees from 1 -1 center at 1 0
\endpicture
\break
{\small\ifx\thesamplenumber\relax % sample?  
Volume \else Sample
\fi\thevolumenumber\ (\thevolumeyear)
\startpage--\finishpage\nl
Published: \publishdate}
\vglue 0.5truein plus 0.4fil minus 0.1truein

{\parskip=0pt\leftskip 0pt plus 1fil\def\\{\par\smallskip}{\ifplaintex\large
\else\Large\fi\bf\thetitle}\par\medskip}   

\vglue 0pt plus 0.1fil 

{\parskip=0pt\leftskip 0pt plus 1fil\def\\{\par}{\sc\theauthors}
\par\medskip}

\vglue 0pt plus 0.1fil 

{\small\parskip=0pt\let\newline\\
{\leftskip 0pt plus 1fil\def\\{\par}{\sl\theaddress}\par}
\expandafter\ifx\theemail\relax    % email address?
\relax\else\vglue 5pt plus 0.02fil minus 2pt\def\\{\stdspace{\rm 
and}\stdspace} 
\cl{Email:\stdspace\tt\theemail}\fi
\ifx\theurl\relax                  % URL given?
\relax\else\vglue 5pt plus 0.02fil minus 2pt\def\\{\stdspace{\rm 
and}\stdspace}
\cl{URL:\stdspace\tt\theurl}\fi\par}

\vglue 7pt plus 0.3fil minus 3pt

{\bf Abstract}
\vglue 5pt plus 0.1fil minus 2pt

\theabstract

\vglue 7pt plus 0.3fil minus 3pt

{\bf AMS Classification numbers}\quad Primary:\quad \theprimaryclass

Secondary:\quad \thesecondaryclass

\vglue 5pt plus 0.3fil minus 2pt

{\bf Keywords}\quad \thekeywords

\vglue 10pt plus 0.5fil minus 5pt

{\small  Proposed: \theproposer\hfill Received: \receiveddate\nl
Seconded: \theseconders\hfill 
\ifx\reviseddate\relax                         % paper revised?
Accepted: \accepteddate                        % no
\else
Revised: \reviseddate                          % yes
\fi}
\eject
}       %  end of definition of \maketitlep

\let\maketitlepage\maketitlep
\let\maketitle\maketitlepage

\font\phead=cmsl9 scaled 950
\font=cmsl9 scaled 1050
\font\pnum=cmbx10 scaled 913
\font\lnum=cmbx10 
\font\pfoot=cmsl9 scaled 950
\font=cmsl9 scaled 1050
\ifplaintex
\headline{\vbox to 0pt{\vskip -4.5mm\line{\small\phead\ifnum
\count0=\startpage ISSN 1364-0380 (on line)
1465-3060 (printed) \hfill {\pnum\folio}\else\ifodd\count0\def\\{ }% 
\ifx\theshorttitle\relax\thetitle\else\theshorttitle\fi\hfill{\pnum\folio}
\else\def\\{ and }{\pnum\folio}\hfill\ifx\theshortauthors\relax\theauthors
\else\theshortauthors\fi\fi\fi}\vss}}
\footline{\vbox to 0pt{\vglue 0mm\line{\small\pfoot\ifnum\count0=\startpage
\copyright\ \gtp\hfill\else
\gt, Volume \thevolumenumber\ (\thevolumeyear)\hfill\fi}\vss
}}
\else
\makeatletter
\def\@oddhead{{\small\lhead\ifnum\count0=\startpage ISSN 1364-0380 (on line)
1465-3060 (printed) \hfill {\lnum\number\count0}\else\ifodd\count0
\def\\{ }\ifx\theshorttitle\relax \thetitle \else\theshorttitle\fi\hfill
{\lnum\number\count0}\else\def\\{ and }{\lnum\number\count0}
\hfill\ifx\theshortauthors\relax 
\theauthors\else\theshortauthors\fi\fi\fi}}\def\@evenhead{\@oddhead}
\def\@oddfoot{\small\lfoot\ifnum\count0=\startpage\copyright\ \gtp\hfill\else
\gt, Volume \thevolumenumber\ (\thevolumeyear)\hfill\fi}
\def\@evenfoot{\@oddfoot}
\makeatother
\fi

\newwrite\gtoutfile
\long\gdef\makeheadfile{  %%% start of definition of \makeheadfile
{\def\\{, }\def\s{ }
\immediate\openout\gtoutfile head.xxx
\immediate\write\gtoutfile{To: math@arxiv.org}
\immediate\write\gtoutfile{Subject: put or rep NNNNN:pppp}
\immediate\write\gtoutfile{--text follows this line--}
\immediate\write\gtoutfile{Proxy-for: \ifx\theasciiauthors\relax
\theauthors\else\theasciiauthors\fi\s<\ifx\theasciiemail\relax\theemail\else\theasciiemail\fi>}
\immediate\write\gtoutfile{\noexpand\\}
\immediate\write\gtoutfile{Authors: \ifx\theasciiauthors\relax
\theauthors\else\theasciiauthors\fi}
\immediate\write\gtoutfile{Title: \ifx\theasciititle\relax
\thetitle\else\theasciititle\fi}
\immediate\write\gtoutfile{Subj-class: GT or SG or MG etc}
\immediate\write\gtoutfile{MSC-class: \theprimaryclass\ifx\thesecondaryclass\relax\else, \thesecondaryclass\fi}
\immediate\write\gtoutfile{Journal-ref: Geom. Topol. \thevolumenumber
(\thevolumeyear) \startpage-\finishpage}
\immediate\write\gtoutfile{Comments: Published by Geometry and Topology at}
\immediate\write\gtoutfile{\s\s http://www.maths.warwick.ac.uk/gt/GTVol\thevolumenumber/paper\thepapernumber.abs.html}
\immediate\write\gtoutfile{\noexpand\\}
\immediate\write\gtoutfile{}
\ifx\theasciiabstract\relax
\immediate\write\gtoutfile{\theabstract}\else
\immediate\write\gtoutfile{\theasciiabstract}\fi
\immediate\write\gtoutfile{}
\immediate\write\gtoutfile{\noexpand\\}
\immediate\write\gtoutfile{}
\immediate\closeout\gtoutfile}}  %%% end of definition of \makeheadfile

\def\maketitlepage{\maketitlep\makeheadfile}
\let\maketitle\maketitlepage

\lognumber{171}
\volumenumber{5}\papernumber{21}\volumeyear{2001}
\pagenumbers{651}{682}
\received{27 March 2001}\revised{17 July 2001}
\accepted{15 August 2001}
\proposed{Yasha Eliashberg}
\seconded{Robion Kirby, Joan Birman}
\published{15 August 2001}

\usepackage{amsmath,amssymb,epsf}

\theoremstyle{plain}
\newtheorem{thm}{Theorem}[section]
\newtheorem{cor}[thm]{Corollary}
\newtheorem{lemma}[thm]{Lemma}
\newtheorem{prop}[thm]{Proposition}

\theoremstyle{remark}
\newtheorem{rem}[thm]{Remark}

 \def\cA{{\mathcal A}}
 \def\cK{{\mathcal K}}
 \def\cC{{\mathcal C}}
 \def\cS{{\mathcal S}}
 \def\cT{{\mathcal T}}
 \def\cV{{\mathcal V}}
 \def\ta{{\mathcal T}_{\subset}}
 
 \def\va{{\mathcal V}_{\subset}}
 \def\bC{{\mathbf C}}
 \def\axis{{\mathbf A}}
 \def\fibr{{\mathbf H}}
 \def\ba{{\mathbf a}}
 \def\bb{{\mathbf b}}
 \def\bc{{\mathbf c}}
 \def\be{{\mathbf e}}
 \def\d{{\delta}} 
  
 \def\e{{\epsilon}} 
  
 \def\L{{\Lambda}}

 \def\s{{\sigma}} 
  
 \def\a{{\alpha}} 
 \def\b{{\beta}} 
  
 \def\r{{\rho}} 
 \def\ra{{\rightarrow}} 
  
 \def\g{{\gamma}}

 \def\z{{\mathbb Z}}
 \def\t{{\tau}}

 \def\pf{\proof}

\begin{document}
\title{On iterated torus knots and transversal knots}
\author{William W Menasco}
\address{University at Buffalo, Buffalo, New York 14214, USA}
\email{menasco@tait.math.buffalo.edu}
\url{http://www.math.buffalo.edu/\char'176menasco}

\begin{abstract}
A knot type is {\em exchange reducible} if an arbitrary closed
$n$--braid representative $K$ of ${\mathcal K}$ can be changed to a
closed braid of minimum braid index $n_{min}({\mathcal K})$ by a
finite sequence of braid isotopies, exchange moves and
$\pm$--destabilizations. (See Figure \ref{figure:destabilizations and
exchange moves}). In the manuscript \cite{[BW]} a transversal knot in
the standard contact structure for $S^3$ is defined to be {\em
transversally simple} if it is characterized up to transversal isotopy
by its topological knot type and its self-linking number.  Theorem 2
of \cite{[BW]} establishes that exchange reducibility implies
transversally simplicity. Theorem \ref{calculus on iterated torus
knots}, the main result in this note, establishes that iterated torus
knots are exchange reducible. It then follows as a Corollary that
iterated torus knots are transversally simple.
\end{abstract}

\asciiabstract{A knot type is exchange reducible if an arbitrary
closed n-braid representative can be changed to a closed braid of
minimum braid index by a finite sequence of braid isotopies, exchange
moves and +/- destabilizations.  In the manuscript [J Birman and NC
Wrinkle, On transversally simple knots, preprint (1999)] a transversal
knot in the standard contact structure for S^3 is defined to be
transversally simple if it is characterized up to transversal isotopy
by its topological knot type and its self-linking number.  Theorem 2
of Birman and Wrinkle [op cit] establishes that exchange reducibility
implies transversally simplicity.  The main result in this note,
establishes that iterated torus knots are exchange reducible.  It then
follows as a Corollary that iterated torus knots are transversally
simple.}

\primaryclass{57M27, 57N16, 57R17} \secondaryclass{37F20}
\keywords{Contact structures, braids, torus knots, cabling, exchange
reducibility}
\maketitlepage
 
\section{Introduction}
\label{section:introduction}
Let $\cC \subset S^3$ be an oriented knot, 
let $\cV_\cC$ be a solid torus neighborhood of $\cC$ and
let $\partial\cV_\cC = \cT_{\cC}
\subset S^3$ be a peripheral torus for $\cC$. The oriented simple closed curve on $\cT_{\cC}$
that represents the homotopy class of $pm+ql$, where $m$ is the meridian homotopy class,
$l$ is the preferred longitude homotopy class and $p,q \in \z$,
is called the {\em $(p,q)$ cable of $\cC$}.  We will use the notion 
$\bC(\cC,(p,q))$ to
indicate the resulting oriented knot of this {\em cabling operation}.
If $\cC$ is the oriented unknot and $(p,q)$ is a co-prime pair then the cabling operation produces a
$(p,q)$--torus knot.  (Our discussion and results can be adapted to the situation where
$(p,q)$ is not co-prime, but since we will be concerned with the iteration of the cabling
operation, to avoid ambiguity we require that the cabling produce a knot.)

We can, of course, iterate the cabling operation.  Starting with
an initial knot $\cC_0$ and a sequence of co-prime $2$--tuples of integers
 $\{(p_1,q_1),(p_2,q_2),\cdots,(p_h,q_h)\}$,
with $p_1 < q_1$,
we can construct the oriented knot
$$\bC(\bC(\cdots\bC(\bC(\cC_0,(p_1,q_1)),(p_2,q_2))\cdots,(p_{h-1},q_{h-1})),(p_h,q_h)).$$
If $\cC_0$ is the oriented unknot then any iteration of the cabling operation produces an
{\em iterated torus knot}.  Letting $(P,Q) = \{(p_1,q_1),(p_2,q_2),\cdots,(p_h,q_h)\}$, the final
iteration produces an oriented knot, $K_{(P,Q)}$, which is on the peripheral torus of
the  next to last knot in the iteration; mainly,
$$\cT_{\bC(\bC(\cdots\bC(\bC(\cC_0,(p_1,q_1)),(p_2,q_2))\cdots),(p_{h-1},q_{h-1}))}.$$

\begin{figure}[ht!]  
\cl{\epsfysize=115pt \epsfbox{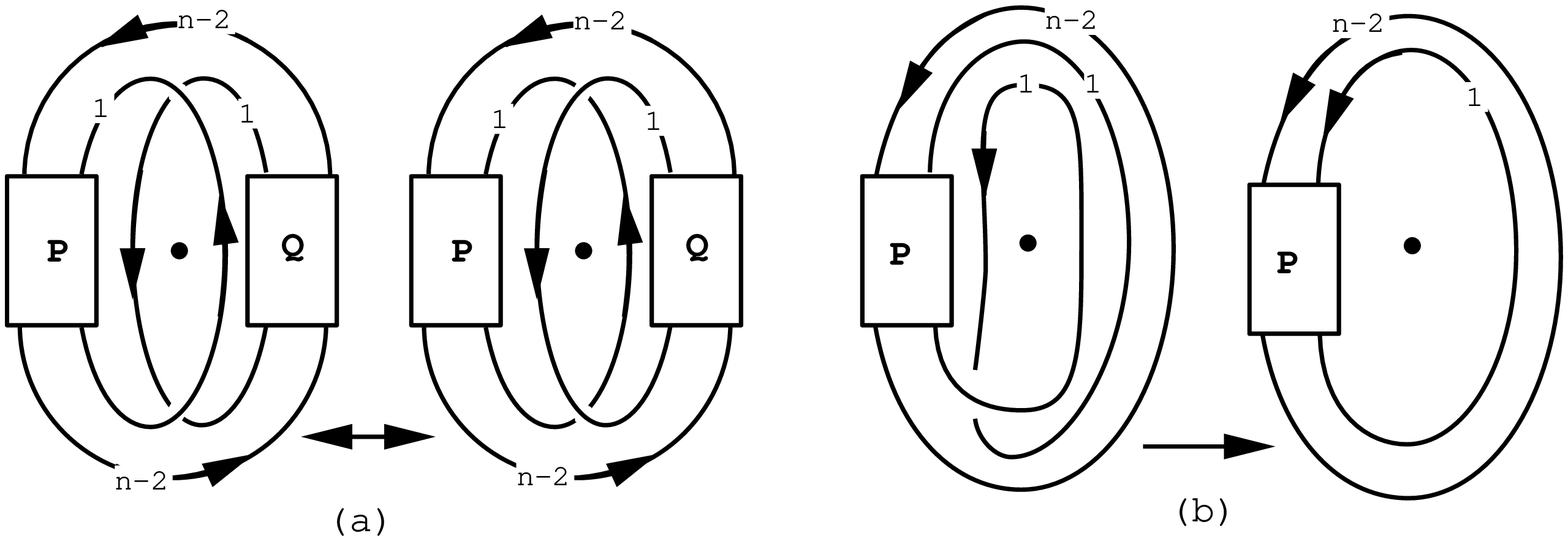}}
\nocolon \caption{}
\label{figure:destabilizations and exchange moves}
\end{figure}

In Section 2 of \cite{[BW]} three moves are discussed which take closed braids to closed
braids,  preserving knot type: braid isotopy, exchange moves and
destabilization.
Braid isotopy means an isotopy in the complement of the braid axis which preserves the braid structure.
The exchange move is a special type of Reidemeister II
move illustrated in Figure
\ref{figure:destabilizations and exchange moves}(a). 
Destabilization means reducing the braid index by eliminating a (positive or
negative) trivial loop, as shown in Figure
\ref{figure:destabilizations and exchange moves}(b). Notice that braid isotopy
and exchange moves preserve both algebraic crossing number and braid
index, whereas
destabilization changes both.
For a more extensive treatment of these isotopies, see
\cite{[BF],[BM2],[BM3]}.)  

As defined in
\cite{[BW]}, a knot type $\cK$ is {\em exchange reducible} if
an arbitrary closed $n$--braid representative $K$ of
$\cK$ can be changed to a closed braid representative of minimum braid index, $n_{min} (\cK)$,
by a finite sequence of braid isotopies, exchange moves and $\pm$--destabilizations.
The main result of \cite{[BM1]} established the exchange reducibility of the unknot.
The main theorem in this paper is an analogous result for iterated torus knots.

\begin{thm} 
\label{calculus on iterated torus knots} 
Iterated torus knots are exchange reducible.
\end{thm}

It should be noted that
in \cite{[BM4]} it was shown that there are knots that can be represented by two non-conjugate
closed $3$--braids where the differing conjugacy classes are not related to each other by an exchange move.
(In fact, they are related to each other by a braid preserving flype.)

The proof of Theorem \ref{calculus on iterated torus knots} involves 
adapting the braid-foliation
machinery developed in \cite{[BM3]} to the situation
where there is a torus in $S^3$ which is being foliated and a knot is embedded
on this torus. It employs the result from
\cite{[Sch]} that an oriented  iterated torus knot $K_{(P,Q)}$ has
an unique braid representative of minimal braid index $\prod_{1}^h p_i$.

Theorem \ref{calculus on iterated torus knots} has an immediate application 
to transversal knots.
Let $\xi$ be the standard contact structure in oriented $S^3$.  The
structure $\xi$ can be thought of as a plane field that is totally
non-integrable.  A knot $K$ is transversal if and only if $K$ intersects each 
plane in the plane field $\xi$ transversally. A transversal isotopy of $K$ is
an isotopy of $K$ in $S^3$ through transversal knots.  (See \cite{[El]}.)  If
$K$ and $K^\prime$ are two transversal knots that are transversally isotopic,
then they are representatives of the same {\em transversal knot type}, $\cT\cK$.

A classical invariant of transversal knot types is a {\em self-linking 
number}, the {\em Bennequin number}, $\beta(\cT\cK)$. The self-linking is defined by
pushing the transversal knot off itself in a direction which is in the contact
plane. A well-defined direction exists because $S^3$ is parallelizable. See
\cite{[BW]} for a  precise description.  A transversal
knot type
$\cT\cK$ is  {\em transversally simple} if it is determined by its topological
knot type $\cK$ and its self-linking number.  In
\cite{[El]} it was first shown that the unknot is transversally simple. In
\cite{[Et]} it was established that positive transversal torus knots are
transversally simple.  In
\cite{[BW]} a more general framework for understanding transversally simple
knots was established.

\begin{thm}[See \cite{[BW]}]
\label{exchange reducible implies transversally simple} 
If $\cT\cK$ is a transversal knot type with associated topological knot type
$\cK$, where $\cK$ is exchange reducible, then $\cT\cK$ is transversally simple.
\end{thm}

Combining Theorems \ref{calculus on iterated torus knots} and
\ref{exchange reducible implies transversally simple}, we have the following immediate corollary. 

\begin{cor}
\label{iterated torus knots are simple} 
Let $\cT\cK_{(P,Q)}$ be a transversal knot type with associated topological
knot type that of the iterated torus knot $K_{(P,Q)}$.  Then $\cT\cK_{(P,Q)}$ is
transversally simple.
\end{cor}

The outline for this note is as follows.
In Section \ref{section:the braid foliation machinery for the torus}
we review and adapt the braid-foliation machinery for the torus that was initially
introduced in \cite{[BM3]}.  We will be concerned with the situation where we are given a torus
which contains a knot $K_{(P,Q)}$ and bounds a solid torus.  However, we do not have a natural way of
identifying the core curve of the solid torus.  Hence, we will use $\ta$ as notation for the given torus
containing $K_{(P,Q)}$.
The foliation machinery on $\ta$ will involve understanding the manipulation of
three different types of foliations---circular, mixed and tiled foliations.  (These foliations will
be defined in Section \ref{section:the braid foliation machinery for the torus}.)
In Section \ref{Proof of theorem 1} we will prove Theorem \ref{calculus on iterated torus knots}
in the special case where
$\ta$ has a circular foliation.
The overriding strategy of the remaining sections
is to reduce the mixed and tiled foliations to circular foliations.
In Section \ref{mixed foliations} we show how destabilizations and exchange moves allow
one to replace a mixed foliation with a circular foliation.  Similarly,
in Section \ref{tiled foliations} we show how destabilizations and exchange moves allow
one to replace a tiled foliation with a circular foliation.

\rk{Acknowledgments} 
This work was partially supported by NSF grant DMS-9626884.
The author wishes to thanks Joan Birman and Nancy Wrinkle
for encouraging him to think about a proof of Theorem \ref{calculus on iterated torus knots}
during his brief sabbatical stay at Columbia University. That stay was
partially supported by NSF grant DMS-9705019 .
 
\section{The braid foliation machinery for the torus}
\label{section:the braid foliation machinery for the torus} 
In this section we adapt the combinatorics of \cite{[BM3]} to the
pair $(K,\ta)$, where $K=K_{(P,Q)}$ is an iterated torus knot which lies on the torus
$\ta$. Since the
exposition in
\cite{[BF]} supplies us with a centralized source for most of the previously
developed machinery, we will use it almost exclusively as our primary reference.
Although the arguments in this note rely heavily on
results in the existing literature, a reader need only consult \cite{[BF]}
and \cite{[BM3]} in almost all cases.
 
Let $K \subset S^3$ be an oriented
closed $n$--braid with axis $\axis$.  Then we can choose
a $2$--disc fibration of the open solid torus $S^3 - \axis $. We will refer to
$ {\bf H} = \{ H_\theta \ | \ \theta \in [0 , 2\pi]\} $ 
as this 2--disc fibration of $S^3 - \axis $. Each $H_\theta$ is a disc with boundary
$\axis$.  We consider the intersection of the $H_\theta$'s with $\ta$---the induced
singular foliation on $\ta$ by $\fibr$.  We have a sequence of lemmas that begin to
standardize this foliation.  These lemmas imitate the similar set of lemmas in
Section 1 of \cite{[BM3]} which dealt with an essential torus in the complement of a 
closed braid.
Since our present case is slightly different (the closed braid is actually a 
homotopically
non-trivial curve on the torus), we will only supply the additional details
needed to adapt the proofs of \cite{[BM3]} to this case.

\begin{lemma}
\label{generic vertices and singular points}
We may assume that:
\begin{itemize}
\item[\rm(i)] The intersections of $\axis$ with $\ta$ are finite in number and transverse.
Also, if $p \in \axis \cap \ta$ then $p$ has a neighborhood on $\ta$ which is radially
foliated by its arcs of intersection with fibers of $\fibr$.
\item[\rm(ii)] All but finitely many fibers $H_\theta \in \fibr$ meet $\ta$ transversally,
and those which do not (the singular fibers) are each tangent to $\ta$ at exactly one
point in the interior of both $\ta$ and $H_\theta$.  Moreover, the tangencies 
(which are contained in singular leaves) are  either local
maxima, or minima, or saddle points.
\item[\rm(iii)] A leaf that does not contain a singular point (a non-singular leaf) is
either an arc having endpoints on $\axis$ or a simple closed curve.
\end{itemize}
\end{lemma}
 
\pf 
We use exactly the same general position argument as in \cite{[BM3]}.\endproof

We refer to the leaves of the foliation of $\ta$ as {\em $\bb$--arcs} and
{\em $\bc$--circles}.
Each $\bb$--arc and each $\bc$--circle lies in both $\ta$ and in some
fiber $H_\theta \in \fibr$.
Since $K \subset \ta$, for all generic $H_\theta \in \fibr$, each point 
of $K \cap H_\theta$
is contained in a $\bb$--arc or $\bc$--circle.  Finally, since $K$
intersects each disc fiber of $\fibr$ coherently, $K$ must intersect each 
non-singular
leaf coherently.

A $\bb$--arc, $b \subset \ta \cap H_\theta$, is {\em essential} if either
$b \cap K \not= \emptyset$,
or both sides of $H_\theta$ split along $b$ are intersected by $K$.
A $\bc$--circle, $c \subset \ta \cap H_\theta$, is {\em essential} if $c \cap K \not= \emptyset$.
The definition of essential $\bb$--arcs and $\bc$--arcs is an adaptation of the definition
in \cite{[BM3]}, however inessential $\bb$--arcs ($\bc$--circles) are still arcs
(respectively, circles) splitting off
sub-discs (respectively, bounding subdiscs) of disc fibers that are not intersected by
$K$.

\begin{lemma}
\label{essential b-arcs and c-circles}
Assume that $\ta$ satisfies (i)--(iii) of Lemma \ref{generic vertices and singular points}.
Then $\ta$ is isotopic to a cabling torus, $\ta^\prime$, such that the foliation of
$\ta^\prime$ also satisfies (i)--(iii) and in addition:
\begin{enumerate}
\item All $\bb$--arcs are essential.
\item All $\bc$--circles are essential.
\item Any $\bc$--circle in the foliation is homotopically non-trivial on $\ta^\prime$.
\end{enumerate}
Moreover, the restriction of the isotopy to $K$ is the identity.
\end{lemma}

\pf
The argument for eliminating inessential $\bb$--arcs is exactly the same as the argument
used in
the proof of Lemma 2 of \cite{[BM3]}.  Similarly, if $c$ is a
$\bc$--circle in the foliation which is homotopically trivial on $\ta$, we must have $c
\cap K = \emptyset$ since $K$ cannot intersect a homotopically trivial circle on
$\ta$ coherently. Then we can, again, employ the argument in the proof of Lemma
2 \cite{[BM3]}, which relies on the fact that $c$ bounds a sub-disc in the
disc fiber that does not intersect $K$.  
\endproof

\begin{figure}[ht!] 
\cl{\epsfxsize=320pt \epsfbox{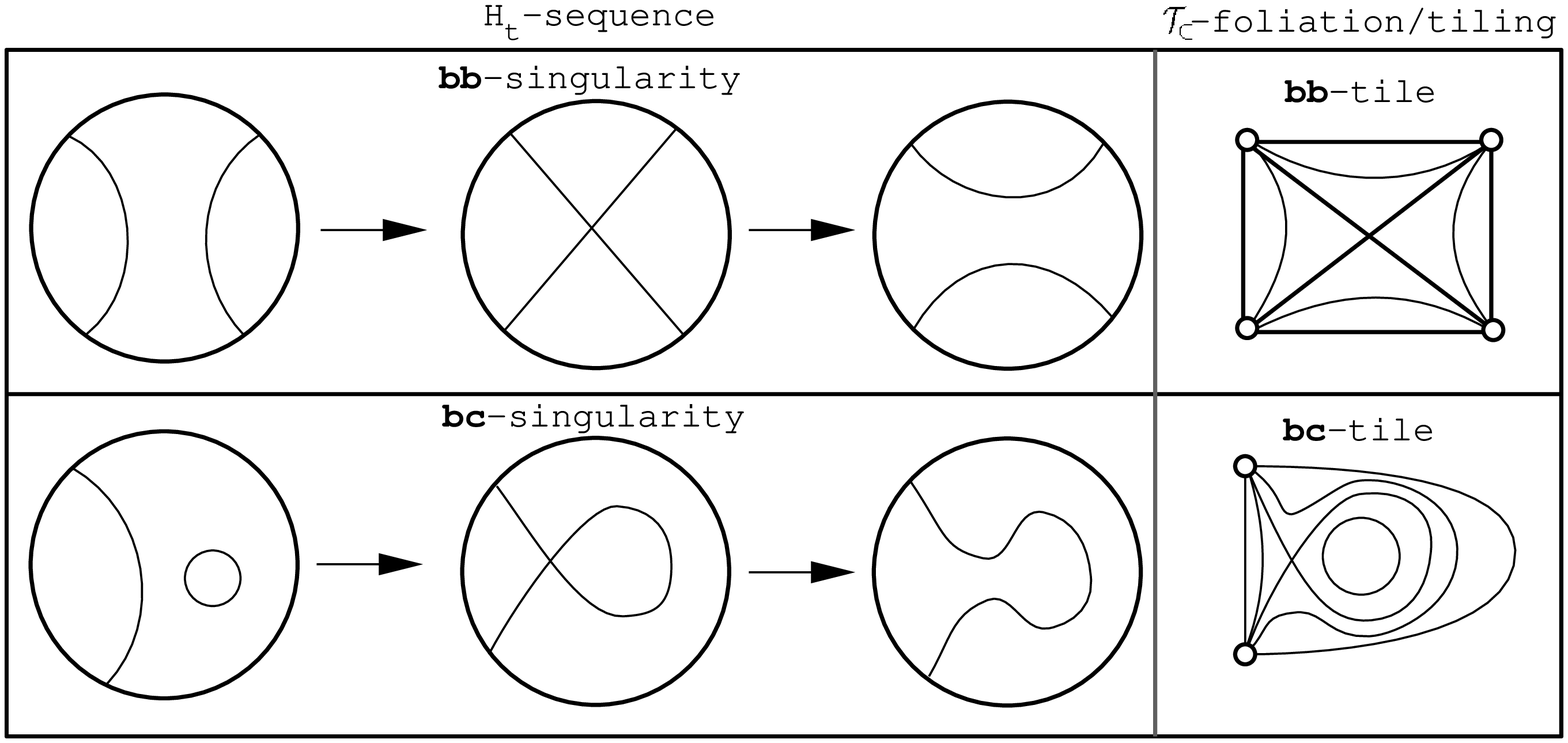}}
\nocolon \caption{}
\label{figure:bb- and bc-tiles}
\end{figure}

We now consider the different types of singularities which can occur in the foliation of $\ta$.
Having two possible non-singular leaves allows for the occurrence of 
three possible
types of singularities:  a $\bb\bb$--singularity resulting from two $\bb$--arcs
``meeting'' at a saddle point;  a $\bb\bc$--singularity resulting from a
$\bb$--arc and a $\bc$--circle forming a saddle point; and $\bc\bc$--singularity
resulting from two $\bc$--circles.  See Figure \ref{figure:bb- and bc-tiles}.

\begin{lemma}
\label{lemma:no cc-singularities}
$\bc\bc$--singularities do not occur.
\end{lemma}

\pf
A surgery between two $\bc$--circles will produce a homotopically trivial $\bc$--circle,
violating statement 3 of Lemma \ref{essential b-arcs and c-circles}.
\endproof

We can now repeat the discussion in Section 1 of \cite{[BM3]} which is based on the
observation that the existence of two different types of singularities implies
the occurrence of two different tiles: $\bb\bb$--tiles and $\bb\bc$--tiles as illustrated
in Figure \ref{figure:bb- and bc-tiles}.  (The reader should note that a $\bb\bc$--tile
is in fact annular.)

\begin{prop}
\label{prop:3 types of foliations}
There are three possibilities for the foliation of $\ta$:
\begin{itemize}
\item A  circular foliation: every leaf is a
$\bc$--circle. There are no singularities.
\item A mixed foliation: there are both $\bb$--arcs and $\bc$--circles. This foliation needs
the occurrence of $\bb\bc$--tiles, but allows the occurrence of $\bb\bb$--tiles.
\item A tiling: that is, a foliation involving only $\bb$--arcs and, thus,
only $\bb\bb$--tiles.  (An example of a tiled torus is given in Figure \ref{figure:example of a tiled torus}.)
\end{itemize}
\end{prop}

At this point an interested reader can reach the punch-line in rapid fashion by reviewing the statements
of Propositions \ref{proposition:mixed implies circular} and \ref{prop:final goal},
and then proceeding to the argument in Section 3.  No understanding
of the machinery in the intervening subsections is needed to understand the proof of
Theorem \ref{calculus on iterated torus knots} in Section 3.  (To understand Remark
\ref{remark:why the proof is difficult} it is
useful to have reviewed the discussion on foliations in Section 2.1 and Section 2.2.)

\subsection{Isotopies of $K$ in $\ta$}
\label{subsection:isotopies of K in ta}
In this section we discuss how the positioning of $K$ in the foliation
of $\ta$ at times allows us to perform destabilizations and exchange moves on
$K$.  The reader should note that any isotopy of $K$ induces an isotopy of $\ta$.

We adopt the terminology of the literature, referring
to the points of $\axis \cap \ta \subset \ta$ as {\em vertices}.  A vertex $v$ is
{\em adjacent} to a leaf in the foliation if $v$ is an endpoint of the leaf.
The {\em valence} of a vertex $v$ is the number of singular leaves which end at $v$.

A {\em meridian} curve $m \subset \ta$ is a curve which intersects $K$ coherently
$p_h$--times and bounds a disc $\Delta_m$ in $S^3$, with $\Delta_m \cap \ta = m$.
We may assume that $m \cap \axis = \emptyset$ and that $m$ is transverse to the leaves
of the foliation on $\ta$.  Being transverse to leaves makes $m$
a closed braid and we orient $m$ so that it has positive linking number with $\axis$.
The orientation on  $m$ is consistent with the forward direction of the
foliation on $\ta$.

The disc $\Delta_m$ is contained in the solid torus $\va$ where $\partial \va = \ta$.
We also assign an orientation
to $\ta$ such that the positive side of $\ta$ ``points'' towards $\va$.
This choice of orientation allows us to make parity assignments to the vertices
and singularities in the foliation of $\ta$.  A vertex or singularity
is {\em positive} if the orientation of $\axis$ agrees with the orientation of the
positive normal vector of $\ta$ at the vertex or singularity, otherwise it is {\em
negative}.  (We will deal with meridian curves and discs in Section 5.  Their introduction
at this point was only necessary for understanding parity assignments to vertices and singularities.)

\begin{figure}[ht!] 
\cl{\epsfysize=100pt \epsfbox{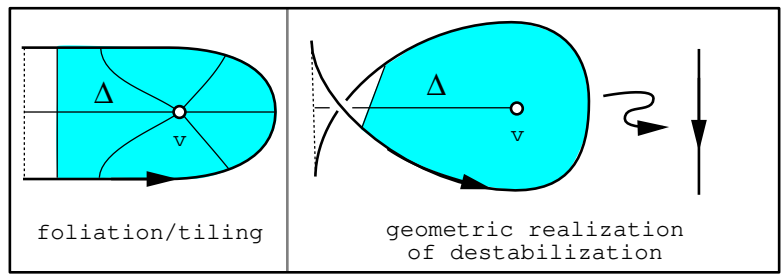}}
\caption{The vertex $v$ has type (a).}
\label{figure:destabilization}
\end{figure}

\begin{lemma}
\label{lemma:destabilization}Let $\Delta \subset \ta$ be a sub-disc such that $\partial \Delta = \g \cup \r $ where
$\r \subset K$ and $\g$ is an arc contained in a singular leaf of the foliation.
Furthermore, suppose that $int(\Delta)$ contains exactly one vertex $v$ and no
singular points.  (See Figure \ref{figure:destabilization}.)  Then an isotopy
of $\r ( \subset K) $ across $\Delta$ to a new position which is transverse to
the leaves of the foliation corresponds to a destabilization of the braid
$K$.
\end{lemma}

\pf
If we split $\ta$ along $K$ we produce an annulus having two
copies of $K$ as its boundary.  The foliation of this annulus will contain a valence
one vertex,
$v$, and the foliated neighborhood of $v$  
will be a {\em type (a) vertex} as
described in Section 2.3 of \cite{[BF]}.  (See Figure \ref{figure:destabilization}.)
As in \cite{[BF]}, the isotopy of $\r$ across $\Delta$ corresponds to a
destabilization of $K$. \endproof

\begin{figure}[ht!] 
\cl{\epsfysize=150pt \epsfbox{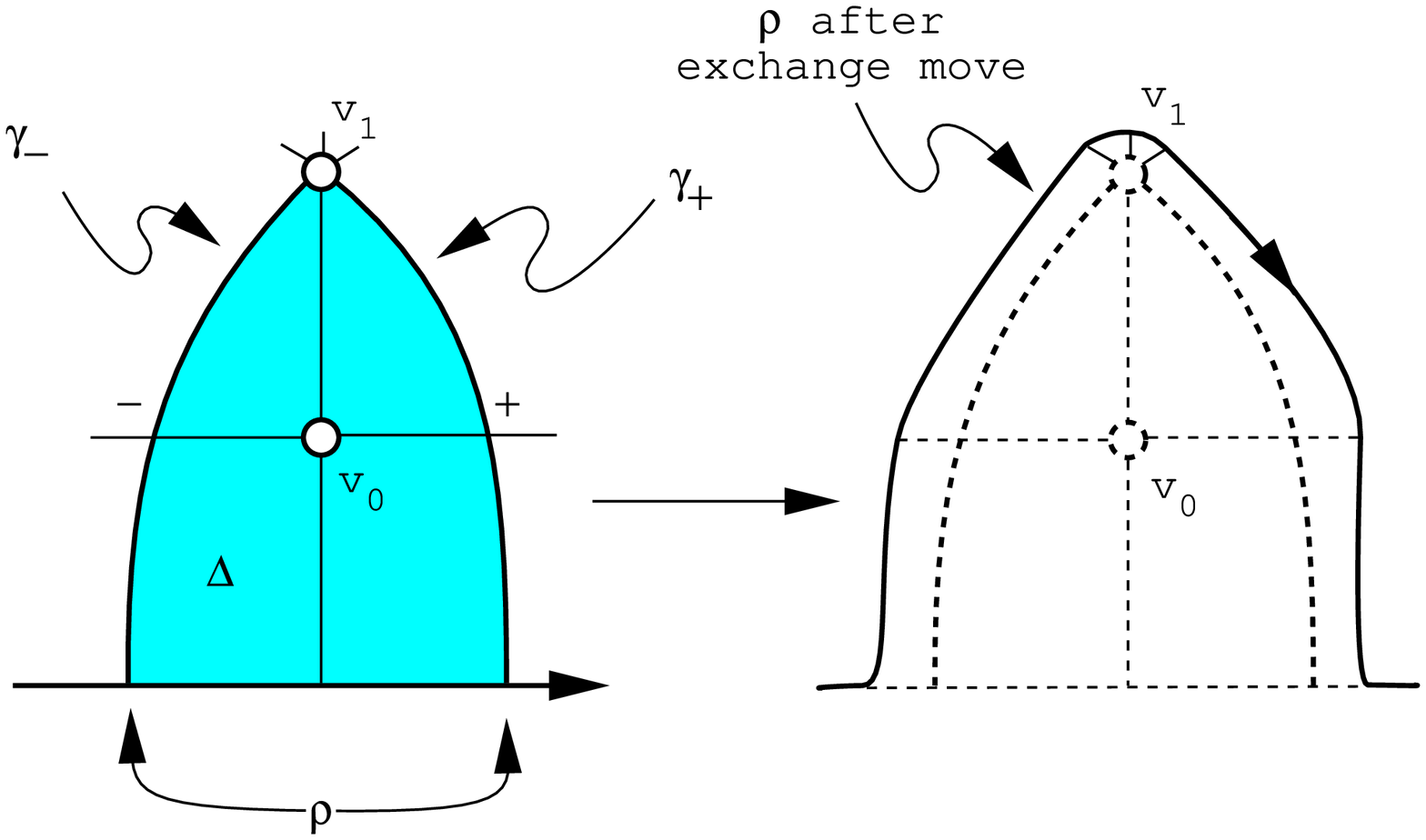}}
\nocolon \caption{}
\label{figure:ab-exchange move}
\end{figure}

\begin{lemma}
\label{lemma:ab-exchange move}
Let $\Delta \subset \ta$ be a sub-disc such that 
$\partial \Delta = \g_+ \cup \g_- \cup \r$,
where  $\g_\pm$ is included in a $\pm$--singular leaf of the foliation and where
$\r$ is a subarc of the knot $K$.  (See Figure \ref{figure:ab-exchange move}.)  Suppose
that the arcs  $\g_+$ and $\g_-$ have one endpoint
at a common vertex $v_1$ and that $int(\Delta)$ contains exactly one vertex
$v_0$ and no singular points.  Then an isotopy
of $\r $ through $\Delta$ to $\r'$ corresponds to an exchange move on the
braid $K$, as illustrated in Figure \ref{figure:destabilizations and
exchange moves}(a).  
\end{lemma}

\pf
Again, we split $\ta$ along $K$ to produce an annulus whose boundary consists
of two copies of $K$.  The foliation of this
annulus will contain a valence two vertex, $v_0$. 
The foliated neighborhood of $v_0$ (see Figure \ref{figure:ab-exchange move})
will correspond to a {\em type (a,b) vertex} as
described in Section 2.4 of \cite{[BF]}.  As in \cite{[BF]}, the isotopy of $\r$ across
such  a neighborhood corresponds to an exchange move. \endproof

\subsection{Manipulating the foliation of $\ta$}
\label{subsection:manipulating the foliation of ta.}
Two operations, {\em change of foliation} and
{\em elimination of a valence two vertex}, played an important role in establishing
control over the foliation of the torus in \cite{[BM3]}.  We now discuss how these
two operations are adapted to our present situation, where we must deal not just
with the foliation of the torus $\ta$, but with the pair $(K,\ta)$.

\begin{figure}[ht!] 
\cl{\epsfysize=185pt \epsfbox{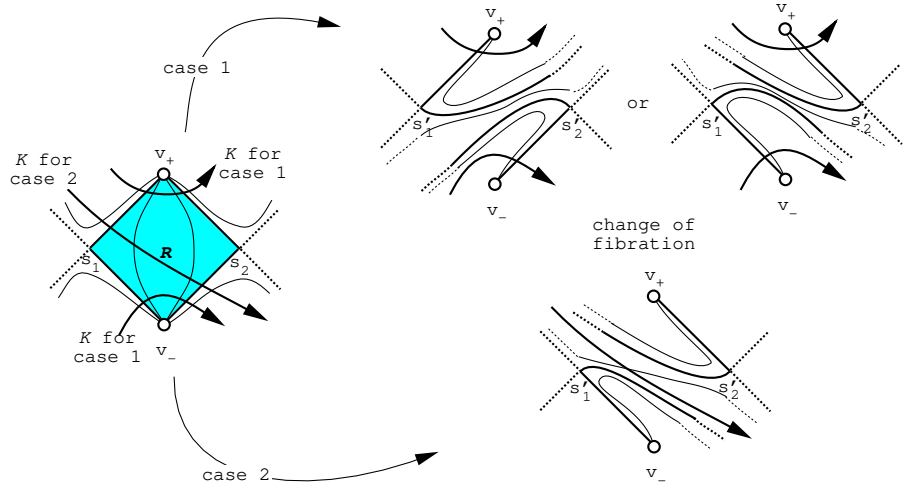}}
\caption{Parity of $s_1$ equals parity of $s_2$.  In case 1 there
are two possible ways to change the fibration of $\fibr$ and, thus, the foliation of
$\ta$.   In case 2 there is only
one way to alter the fibration.}
\label{figure:change of fibration}
\end{figure}

\noindent
{\bf Change of foliation}\qua Let $R \subset \ta$
be the topological disc that is closure of a 
connected region foliated by $\bb$--arcs having common endpoints at vertices
$v_+ , v_- \subset \ta$.  (The subscripts indicate the parity.)
Let $s_1 , s_2 \subset R \subset \ta$ be the two singularities
that are on the boundary of $R$ and assume that their parity is the same.
The existence of such a region $R$ is the central assumption in the discussion in
Section 2.1 of \cite{[BF]}.  Specifically, Theorem 2.1 of \cite{[BF]} allows the two
singular points to be either $\bb\bc$ or $\bb\bb$ singularities.  Figure
\ref{figure:change of fibration} shows how the application of this result from \cite{[BF]}
alters the foliation of $\ta$.
To adapt the \cite{[BF]} change of foliation to our present situation 
we need only check
that the presence of the knot $K \subset \ta$ does not obstruct the 
change in foliation.

How might $K$ intersect $R$?  We consider an arc $\alpha \subset K \cap R$.
Since $K$ is transverse to the foliation
of $\ta$ we know that $\alpha$ is transverse to all of the $\bb$--arcs of $R$ and
is away from the singular and vertex points of $R$.
We can list the possible cases for a component $\alpha \subset K \cap R$ as follows:

\begin{enumerate}
\item The arc $\alpha$ splits the disc $R$ into two discs, one containing
the two singular points and a vertex, and the other containing only one vertex.
\item The arc $\alpha$ splits the disc $R$ into two discs, each of which contains
a vertex and a singular point.
\end{enumerate}

We recall (see \cite{[BF]}) that there are two possible ways the foliation
of
$R$ can be altered. If case 1 occurs then either of these foliation changes is
permissible.  (See Figure
\ref{figure:change of fibration}.)  If case 2 occurs then only one of the changes in
foliation is possible because only one of the changes results in $K$ still being
transverse to the foliation.  Figure
\ref{figure:change of fibration} illustrates how a case 2 arc $\alpha$ determines
the local change in foliation in $R$. 
The proof that these changes in foliation correspond to braid isotopies of the
$(K,\ta)$ pair is straight forward, but has numerous details.  We will not
repeat the argument in \cite{[BF]}.  We will refer to this braid isotopy (which
only alters the foliation of $\ta$ in a disc neighborhood of $R$) as a {\em
change of foliation}. The following lemma
describes the main features of the change in foliation that we need:

\begin{lemma}
\label{lemma:change of fibration}
Let $R$ be the closure of a region foliated by $\bb$--arcs such that
there are vertices $v_+,v_- \subset \partial R$ and singularities 
$s_1, s_2 \subset \partial R$.
Then there exists a change of foliation such that the valence of these two vertices
has been decreased.
\end{lemma}

\begin{figure}[ht!]
\cl{\epsfysize=130pt \epsfbox{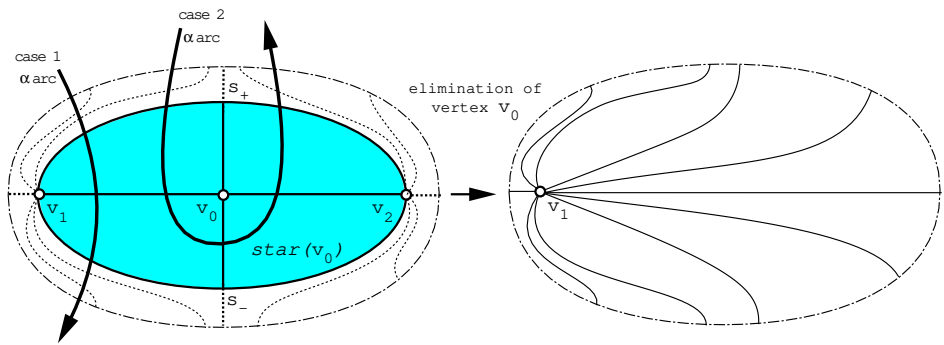}}
\caption{The vertex $v_0$ is a type (b,b) vertex.
In the situation where the case 2 $\a$ arcs are absent, after an exchange move,
the $\bb$--arcs between $v_0$ and $v_2$ are inessential.}
\label{figure:valence two vertex}
\end{figure}

\noindent
{\bf Elimination of valence two vertices}\qua
We next consider the  configuration in the left sketch in 
Figure \ref{figure:valence two vertex}.  (For the moment the reader should ignore the arcs labeled
$\a$ in Figure \ref{figure:valence two vertex}.  They will be referred to in the proof of Lemma
\ref{lemma:eliminating valence two vertices}.)

Let $v_0 \subset \ta$  be a vertex of valence $2$.  Let $star(v_0)$ be
the topological disc which is the closure of the union of all the 
$\bb$--arcs having $v_0$ as an endpoint.  Let $\e$ be the parity of $v_0$ then 
$star(v_0)$ contains exactly two additional vertices, $v_1$ and $v_2$, which will
necessarily have parity $-\e$.  By Lemma 3.1 of
\cite{[BF]}, $star(v_0)$ must contain both $+$ and $-$ singular points.
So let $s_+, s_- \subset star(v_0)$ be the singular points of positive and 
negative parity, respectively.  In this situation we have the following lemma.

\begin{lemma}
\label{lemma:eliminating valence two vertices}
Let $v_0$ be a valence two vertex, with $star(v_0)$ topologically a disc.
Suppose $star(v_0)$ contains vertices $v_1\cup v_2$ and singularities
$s_+\cup s_-$. Then, after an isotopy of $K$ involving only exchange moves, 
destabilizations and braid
isotopies, we can locally eliminate vertices $v_0$ and $v_1$ (or $v_0$ and $v_2$) and
singularities $s_+$ and $s_-$ in the foliation of $\ta$, as illustrated in the
right sketch in Figure \ref{figure:valence two vertex}.  Moreover, away from
$star(v_0)$ the foliation of $\ta$ remains unchanged.
\end{lemma} 

\pf By Theorem 2.2 of \cite{[BF]}, we know that if $K \cap star(v_0) =
\emptyset$ then it is possible to perform an exchange move on $K$ such that the
foliation of $\ta$ remains unchanged (except for a change in the cyclic ordering of
the vertices along $\axis$). After the exchange move, the $\bb$--arcs of $star(v_0)$
will be outermost, ie, for $\bb$--arc $b \subset star(v_0)$, if $ b \subset
H_\theta \in \fibr $ then
$b$ splits off a sub-disc $\Delta \subset H_\theta$
such that $int(\Delta) \cap K = \emptyset$.  Furthermore, if $K \cap star(v_0) =
\emptyset$ then
$\ta$ can be isotopied along $\Delta$, eliminating $v_0$ and, say, $v_2$, and the
two singularities, $s_+$ and $s_-$.  The foliation away from $star(v_0)$ remains 
unchanged, but globally it has two fewer vertices and two fewer singularities.

How might we eliminate $v_0$ if $K \cap star(v_0) \not= \emptyset$?
As with the change of foliation, we consider an arc $\alpha \subset K \cap star(v_0)$.
Since $K$ is transverse to the foliation
of $\ta$ we know that $\alpha$ is transverse to all of the $\bb$--arcs in $star(v_0)$
and does not intersect the singular points and vertices of $star(v_0)$.
We list the possibilities for a component $\alpha \subset K \cap star(v_0)$:
\begin{enumerate}
\item The arc $\alpha$ splits the disc $star(v_0)$ into two discs, one containing
the two singularities and two of the three vertices, and the other containing 
only one vertex.  (See Figure \ref{figure:valence two vertex}.)
\item The arc $\alpha$ splits the disc $star(v_0)$ into two discs, one containing
$v_0$ and a single singularity (say $s_+$), and the other containing 
$v_1$, $v_2$ and $s_-$.  (See Figure \ref{figure:valence two vertex}.)
\end{enumerate}

To deal with case 1, we must in fact consider two situations:  there exists $\bb$--arcs
of $star(v_0)$ that do not intersect $K$; and $K$ intersects every $\bb$--arc of $star(v_0)$.
If some $\bb$--arc of $star(v_0)$ is not intersected by $K$ then after the exchange move
there will necessarily be $\bb$--arcs that are inessential.  The situation where
$K$ intersects every $\bb$--arc of $star(v_0)$ requires a little more work.

Notice that if $K$ intersects every $\bb$--arc of $star(v_0)$, there is necessarily an arc
$\g \subset \partial(star(v_0))$ contained in a singular leaf and having endpoints $v_1 ,
v_2$ which $K$ must intersect incoherently.  The topology of $\ta$, thus, forces
the existence of a sub-arc $\g^\prime \subset \g$, a path $\r \subset K$, 
and a sub-disc
$\Delta \subset \ta$ such that:  $\Delta \cap \g = \g^\prime$; and 
$\Delta \cap K = \r$.
Since $\Delta \subset \ta$, $\Delta$ inherits a foliation that basically mimics
that of
a Seifert disc for the unknot.  Theorem 4.3 of \cite{[BF]} allows us to isotop
$\r$ through $\Delta$ using exchange moves and destabilizations until $K$ is moved
off $\g$ and we have an arc $\alpha \subset star(v_0)$ that corresponds to
case 2.

Finally, to deal with case 2, notice that the sub-disc containing $v_0$ that $\alpha$
splits off has a foliation corresponding to that of a valence one vertex.  
(See Figure \ref{figure:valence two vertex}.)  This is the configuration in
Lemma \ref{lemma:destabilization}.  We can, thus, destabilize $K$ to remove
a case 2 intersection arc.
\endproof

\section{Proof of Theorem \ref{calculus on iterated torus knots} in the case of a circular foliation on $\ta$}
\label{Proof of theorem 1}
In this section we assume that $K$ is contained in a torus $\ta$ that does not intersect the axis
$\axis$, ie, $\ta$ has a circular foliation. With this assumption
our argument for proving Theorem \ref{calculus on iterated torus knots}
is inductive and we need the following result.

\begin{cor}
\label{corollary:cable of exchange reducible is exchange recucible}
Let $\cK$ be an exchange reducible knot type.  Then $\bC(\cK,(p,q))$ is also exchange reducible.
\end{cor}

\pf
Let $K_{(p,q)}$ be any closed braid representative of $\bC(\cK,(p,q))$ and 
$\ta$ be a cabling
torus.  By Proposition \ref{prop:final goal}
we know that if the foliation of $\ta$ is a tiling then through
a sequence of exchange moves and destabilizations
we can replace this tiling with a mixed foliation.
By Proposition \ref{proposition:mixed implies circular}
we can, through a sequence of exchange moves and destabilizations,
replace a mixed foliation of $\ta$ with a circular foliation.  The core of 
this circularly
foliated torus is a braid representative of the cabling knot $\cK$ which we call $K$.
By assumption $K$ represents an exchange reducible knot type.

Now, referring back to Figure \ref{figure:destabilizations and exchange moves},
we notice that we can alter the destabilizing move in (b)
by replacing the weight of $1$ with a weight of $p$, ie, we think of $p$ parallel strands
instead of $1$ strand.  Similarly, in (a)
we can replace the weight of $1$ on the strands involved in the exchange move 
isotopy with a 
weight of $p$.  Thus, a destabilization of $K$ results in $p$ destabilizations of $K_{(p,q)}$,
and an exchange move on $K$ results in an exchange move on $K_{(p,q)}$---the peripheral torus
of $K$ is the cabling torus circularly foliated.  By a classical result in \cite{[Sch]} (cf Satz 23.2),
once $K$ is of
minimal braid index, we will have $K_{(p,q)}$ achieving its minimal braid index when the cabling torus
is circularly foliated.
\endproof

Now, let $U$ be the unknot.  By Theorem 1 of \cite{[BM1]}, $U$ is exchange reducible.
Then, by Corollary \ref{corollary:cable of exchange reducible is exchange recucible}
$\bC(U,(p_1,q_1))$ is exchange reducible.

Inductively, suppose that $$\bC(\cdots\bC(\bC(\cC_0,(p_1,q_1)),(p_2,q_2))\cdots,(p_{i-1},q_{i-1}))$$
is exchange reducible.  Then, again, by Theorem 1 of \cite{[BM1]}
$$\bC(\bC(\cdots\bC(\bC(\cC_0,(p_1,q_1)),(p_2,q_2))\cdots,(p_{i-1},q_{i-1})),(p_i,q_i) ) $$
is exchange reducible.
Thus, $K_{(P,Q)}$ is exchange reducible and Theorem \ref{calculus on iterated torus knots} is
established for this special case of the foliation of $\ta$.

\begin{figure}[ht!]
\cl{\epsfysize=320pt \epsfbox{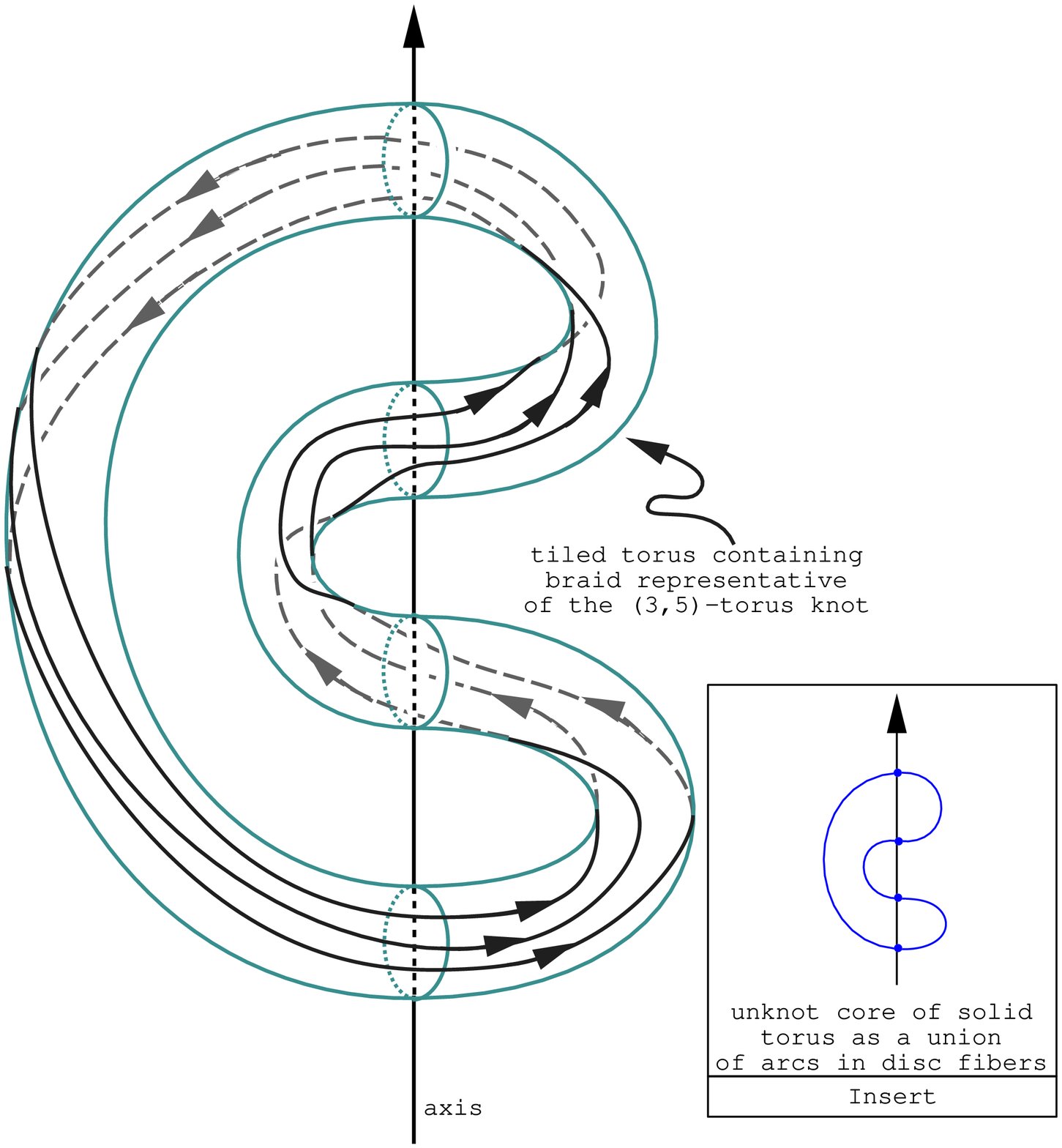}}
\nocolon \caption{}
\label{figure:example of a tiled torus}
\end{figure}

\begin{rem}
\label{remark:why the proof is difficult}
{\rm The ease with which we are able to prove Theorem \ref{calculus on iterated torus knots}
for this special case is due to the fact that
the core of $\va$ can be canonically chosen as a closed braid in the complement of $\axis$, ie,
$\cC \subset \va$ can be taken as the union of the ``centers'' of the
sub-discs in the disc fibers of $\fibr$ that are bounded by $\bc$--circle.  The difficulty with the remaining
two cases---mixed and tiled foliations---is that there is no similar canonical choice for the core of $\va$.
In Figure \ref{figure:example of a tiled torus}, we give an example of a tiled $\ta$ containing the $(3,5)$--torus
knot $K_{(3,5)}$.  This $\ta$ can be seen as the peripheral torus of an
unknot that is represented by the union of four arcs that have their endpoints on $\axis$ and are
contained in disc fibers of $\fibr$.  The insert in Figure \ref{figure:example of a tiled torus} depicts
this representation of the unknot---a core of $\va$.  Since this core intersects $\axis$ at four points,
$\ta$ will necessarily intersect $\axis$ in eight points, ie four times two.  Because of Euler characteristic
considerations, the tiling of $\ta$ will then have eight $\bb$ tiles. 
Notice that there are several possible ways we can $\e$--push this unknotted core off $\axis$ so that
resulting $\cC$ is transverse to $\fibr$.  In particular, it is possible to produce a transverse $\cC$ such
that any $\bb$--arc $ b \subset \ta \cap H_t \subset H_t$ with $ b \cap K_{(3,5)} = \emptyset$
splits the set $\cC \cap H_t$ in $H_t$, for $H_t \in \fibr$.  That is, 
all of the $\bb$--arcs of $\ta$ are {\em essential} in the complement of $\cC \sqcup K_{(3,5)}$.  Using
an exchange move or destabilization isotopy of $\cC$ to induce an isotopy of $\ta$ then
becomes unworkable.  Thus, the key
strategy in the remaining two cases will be to find inessential $\bb$--arcs on $\ta$ in the knot complement
that will enable us to simplify $\ta$.}
\end{rem}

\section{Replacing mixed foliations with circular foliations}
\label{mixed foliations}
In this section we assume that we are given the braid--torus pair $(K,\ta)$, and that
the induced foliation on $\ta$ satisfies the conclusions of Lemmas
\ref{generic vertices and singular points},
\ref{essential b-arcs and c-circles} and \ref{lemma:no cc-singularities}.
Moreover, we assume that the foliation of $\ta$ contains $\bc$--circles. Our goal is to
reduce this case to the special case when the foliation contains only $\bc$--circles.

\begin{prop}
\label{proposition:mixed implies circular}
Let $\ta$ have a mixed foliation.  Then after a sequence of exchange moves 
and destabilizations of $K$ and isotopies of $\ta$ the foliation of $\ta$ may be
assumed to be circular.
\end{prop}

\pf
Following the discussion in Section 3 of \cite{[BM3]}, we observe that 
a singularity between a $\bb$--arc and a $\bc$--circle foliates a 
{\em$\bb\bc$--annulus}.
Moreover, since a $\bc$--circle must have a $\bb\bc$--singularity
both in the forward and backward direction, these $\bb\bc$--annuli 
occur in pairs.

\begin{figure}[ht!] 
\cl{\epsfysize=150pt \epsfbox{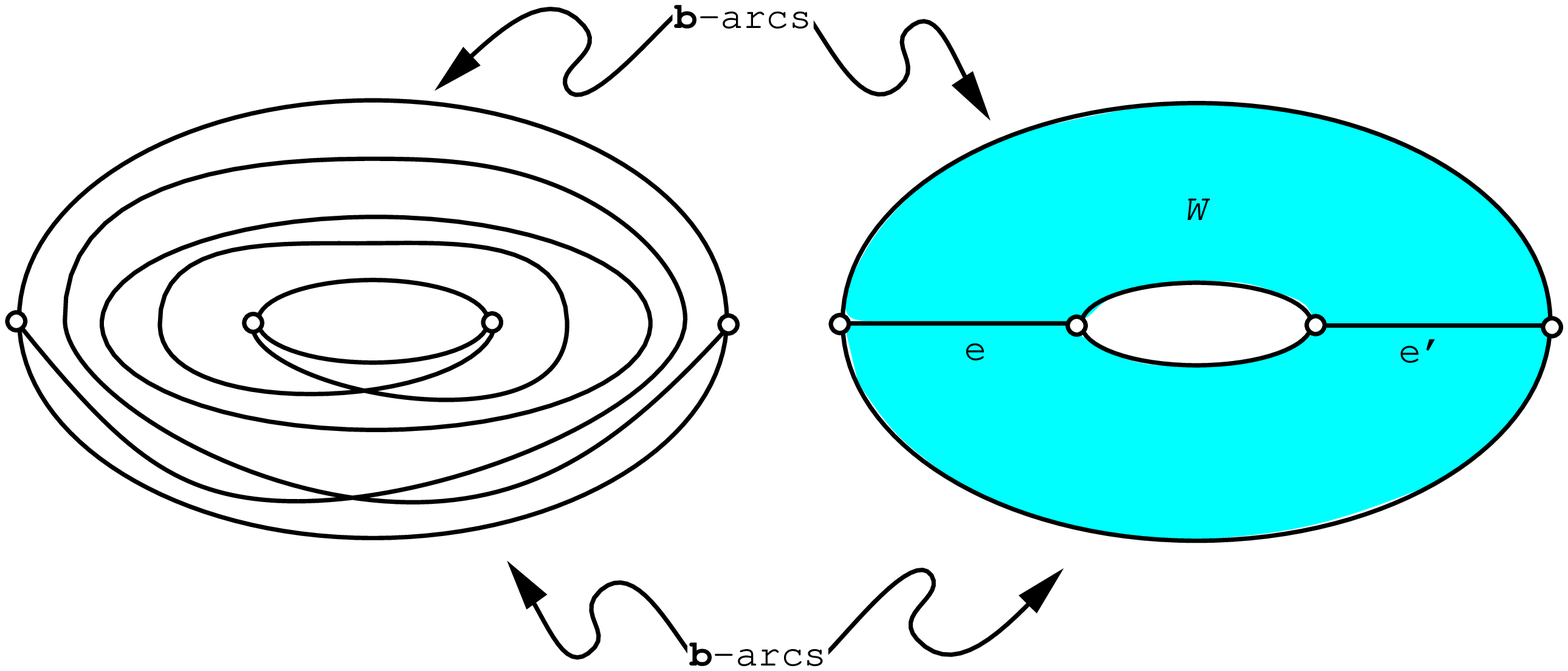}}
\nocolon \caption{}
\label{figure:bebe tiles}
\end{figure}

If we adjoin two $\bb\bc$--annuli
along their common $\bc$--circle we will produce an annular region, $W$, that has each boundary
curve the union of two $\bb$--arcs and, thus, contains two vertices.
See Figure \ref{figure:bebe tiles}.

We can cut $W$ open along two new disjoint edges,
$e$ and $e^\prime$, each having its endpoints at vertices which are on distinct
components of $\partial W$, as in the bottom of Figure \ref{figure:bebe tiles}.
We see that $W$ is the union of two $\bb\be$--tiles where the boundary of a $\bb\be$--tile
has two $\bb$--arcs and two of the new $\be$--arcs.  The number of $\bb\be$--tiles
constructed in a mixed foliation is exactly equal to the number of
$\bb\bc$--annuli.  The vertices in the $\bb\be$--tiling still correspond to the points
of $\axis \cap \ta$.

The cellular decomposition of $\ta$ coming from the $\bb\be$--tiling yields an 
``Euler characteristic'' formula (equation (1) of \cite{[BM3]}).  Namely,
let $V(\b,\e)$ be the number of vertices in the $\bb\be$--tiling of $\ta$ that
is adjacent to $\b$ $\bb$--arcs and $\e$ $\be$--arcs.
Then we have:
$$ 2 V(2,0) + 2V(1,1) + V(2,1) + V(3,0) = \sum_{i=5}^\infty \sum_{\b = [{i}/{2}]}^i
(i-4) V(\b , i - \b) \eqno (1) $$
where both the left hand side and the right hand side are non-negative.

Referring to the discussion in \cite{[BM3]} we recall that $V(1,1)=V(2,1)=0$,
since both a $\bb\be$--vertex (respectively, $\bb\bb\be$--vertex) cannot be geometrically
realized. A $\bb\be$--vertex (respectively, $\bb\bb\be$--vertex) is one that is cyclically
adjacent to
$\bb$--arcs and then $\be$--arcs (respectively, $\bb$--arcs and then, after a singularity, to
$\bb$--arcs and then
$\be$--arcs ). This vertex notation generalizes in the obvious manner.
If $V(3,0)\not=0$ then we have a $\bb\bb\bb$--vertex.  Such a vertex, $v$, must be
adjacent to singular leaves of both positive and negative parity
(cf Lemma 3.1, \cite{[BF]}).  Thus, since $v$ is of odd valence, $star(v)$ will contain a sub-disc region
satisfying the assumptions of Lemma \ref{lemma:change of fibration}.  We can
then perform the change of foliation decreasing the valence of $v$ so that it
becomes a $\bb\bb$--vertex.

If $V(2,0)\not=0$ then there exists a vertex, $v$, that is a $\bb\bb$--vertex.
$star(v)$ will satisfy the assumption of Lemma \ref{lemma:eliminating valence two
vertices}. After possibly some sequence of exchange moves and destabilizations, we can
simplify the tiling of $\ta$ so that $V(2,0) = 0$.

\begin{figure}[ht!] 
\cl{\epsfysize=120pt \epsfbox{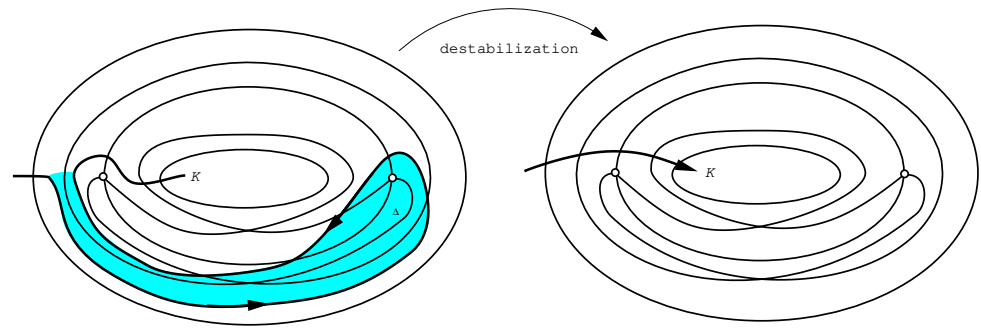}}
\caption{$K$ enters into the foliated region and intersects
a singular leaf twice.  This allows for a destabilization of $K$.  
The shaded region $\Delta$
corresponds to the $\Delta$ of Lemma \ref{lemma:destabilization}.}
\label{figure:bebe vertice destabilization}
\end{figure}

Since we can now assume the left side of equation (1) is zero, we can assume that
the only possible vertices in the tiling of $\ta$ are $\bb\be\bb\be$--vertices.
Figure \ref{figure:bebe vertice destabilization} illustrates the foliation of the
annular neighborhood around such valence four vertices.  Notice that
the $\bc$--circles must be intersected coherently $p_h$--times by $K$.
If $K$ does intersect all of the $\bb$--arcs of this local foliation then
there will exist a sub-disc satisfying the assumptions of Lemma
\ref{lemma:destabilization}.
We can then destabilize $K$ so that, as in Figure
\ref{figure:bebe vertice destabilization}, there will be $\bb$--arcs that $K$ does
not intersect.  An innermost such annulus (one
having consecutive vertices on $\axis$)  will then have an inessential $\bb$--arc.
We may eliminate it by an isotopy of $\ta$.
\endproof

\section{Replacing tilings with mixed foliations}
\label{tiled foliations}
As in the previous section, we start by assuming that we are given
the braid--torus pair $(K,\ta)$ such that
the induced foliation on $\ta$ satisfies the conclusions of Lemmas
\ref{generic vertices and singular points}, \ref{essential b-arcs and c-circles} and \ref{lemma:no cc-singularities}.
But, in this section, we assume that the foliation of $\ta$ contains no $\bc$--circles.
Our goal in this section is to prove that after a sequence of exchange moves,
braid isotopies and  destabilizations of $K$, and isotopies of $\ta$, we may
assume that the foliation of $\ta$ can be assumed to be mixed.

Recall $m \subset \ta$ is a meridian curve and $\Delta_m \subset \va$ is a meridian disc.
Since $\Delta_m$ is a spanning surface having a closed braid boundary, the induced
foliation on the disc $\Delta_m$ satisfies the conclusion of Theorem 1.1, part (i), of
\cite{[BF]}.  We review the most important features of this foliation, all of which are
developed in detail in \cite{[BF]}. If we consider $\Delta_m \cap H_\theta \subset
H_\theta$, where
$H_\theta$ is a generic disc fiber, we see that this intersection
contains two types of arcs:  the first, called $\ba$--arcs,  have one endpoint on
$\axis$ and one endpoint on $m$; the second, called
$\bb$--arcs, have both endpoints on $\axis$.  The singularities
are of three types:
$\ba\ba$--singularities (saddle-point singularities formed by two $\ba$--arcs);
$\ba\bb$--singularities (saddle-point singularities formed by an $\ba$--arc and a
$\bb$--arc); and $\bb\bb$--singularities (saddle-point singularities formed by
two $\bb$--arcs).  Note that the orientation on $m$ induces an orientation on
$\Delta_m$, and that this allows us to assign parities to the vertices (points in
$\Delta_m\cap\axis$) and singularities in the foliation in the same
manner as we did in Section \ref{section:the braid foliation machinery for the torus}
for $\ta$.

The argument in this section will involve ``simplifying'', first, the
pair $(K,\ta)$ then, second, the triple $(K,\ta,m)$ and, finally, the quadruple
$(K,\ta,m,\Delta_m)$.  Our goal is to prove the following proposition.

\begin{prop} 
\label{prop:final goal}
Let $(K,\ta)$ be a knot--torus pair such that the
foliation of $\ta$ is a tiling. Then, after a sequence of exchange moves and
destabilizations of $K$, the foliation of $\ta$ may be assumed to be mixed.
\end{prop}

\subsection{The knot--torus pair}
\label{subsection:the knot-torus pair}
In this subsection we define the {\em complexity of $(K,\ta)$} to be
$\L(K ,\ta)\! =\! (n_1, n_2)$, where
$n_1$ is the braid index of $K$ and 
$n_2$ is the number of vertices in the foliation of $\ta$.
We use lexicographical ordering on the $2$--tuples $(n_1,n_2)$ to
give an ordering on $\L(K , \ta)$.

Following the discussion in Section 2 of \cite{[BM3]}, we recall that the
foliation of $\ta$ (absent of $\bc$--circles) yields a tiling of $\ta$ by
$\bb\bb$--tiles.  This tiling in turn gives us a 
cellular decomposition of $\ta$ and, thus, an 
``Euler characteristic'' formula (compare with equation (1) of \cite{[BM3]}).  
Let $V(\b)$ be the number of vertices in the $\bb\bb$--tiling of $\ta$ that
are adjacent to $\b$ $\bb$--arcs.  
We have:
$$ 2 V(2) + V(3) = V(5) + 2V(6) + 3V(7) + \cdots  \eqno (2) $$
where both the left hand side and the right hand side are non-negative.

If $V(3)\not= 0$ then we have (as in \cite{[BM3]}) a $\bb\bb\bb$--vertex, $v$.  Since
$star(v)$ must contain both positive and negative singularities
(see Lemma 3.1 of \cite{[BF]}), it must contain a sub-disc region
satisfying the assumptions of Lemma \ref{lemma:change of fibration}.
After a change of foliation we can assume that the valence of $v$ is two. We can
then apply the procedure in Lemma \ref{lemma:eliminating valence two vertices} to
reduce $n_2$ and, thus, $\L(K,\ta)$.

If $V(2) \not= 0$ then we can apply the procedure in
Lemma \ref{lemma:eliminating valence two vertices} straight away to reduce 
$\L(K,\ta)$.  After some number of changes of these two types we may assume that the
LHS of Equation (2) is zero. But then the RHS is too. We can conclude that the only
possibility is that   $V(4)$
is non-zero.  Moreover, if $v \subset \ta$ is a $\bb\bb\bb\bb$--vertex then the
four singular leaves which intersect at $v$ are $+,-,+,-$ in that cyclic order.
(Otherwise, a change of foliation could be performed to reduce the valence of
$v$ to three.)  The salient features of this {\em standard tiling} are:  all
vertices of $\ta$ are valence four; the parity pattern on the tiling is a
``checkerboard'' pattern; and all $\bb$--arcs are essential.  The following 
proposition summarizes the above
discussion.

\begin{prop}
\label{proposition:ta has a standard tilng}
Let $(K,\ta)$ be a knot--torus pair where the foliation of $\ta$ is a tiling.
If the tiling of $\ta$ is non-standard then, through the use of exchange moves,
we may replace the knot--torus pair with $(K^\prime, \ta^{\prime})$ such that
$\ta^{\prime}$ has a foliation that is either a standard tiling or a mixed foliation.  
Moreover, $\L(K^\prime,\ta^{\prime}) < \L(K,\ta)$.
\end{prop}

Next, we define four graphs, $G_{+,+}$, $G_{+,-}$, $G_{-,+}$ and $G_{-,-}$ in $\ta$.
The vertices of $G_{\d,\e}$ are the vertices of $\ta$ having parity $\d$ and the edges of
$G_{\d,\e}$ are sub-arcs of singular leaves which join the two $\d$ vertices
in the $\bb\bb$--tiles having a singularity of parity $\e$.  Notice that
the definition of these graphs forces $G_{+,+}$ ($G_{+,-}$) to be disjoint
from $G_{-,-}$ (respectively, $G_{-,+}$).  Moreover, the parity checkerboard pattern
to the standard tiling implies that each vertex of any of the four graphs is
adjacent to exactly two graph edges.  Thus, each component of $G_{\d,\e}$ is a
simple closed curve (scc)
on $\ta$ which, by Lemma 3.8(i) of \cite{[BF]}, is also
homotopically non-trivial on $\ta$.  (For a more complete analysis of $G_{\d,\e}$,
see Section 3 of \cite{[BF]}.)

\begin{lemma}%%%%\hbox
\label{lemma:graph implies destabilization and exchange moves}
Let $\cC \subset \ta$ be a curve such that $\cC \in \{K,m\}$.
Then we can apply either a destabilization or an exchange move to $\cC$ in the
following situations.
\begin{itemize}
\item[\rm(a)] Suppose that as $\cC$ is traversed, a sub-arc $\r \subset \cC$ has the
following sequential intersections with our four graphs: $G_{\d,\e} \ra G_{-\d,-\e}
\ra G_{\d,-\e} \ra G_{-\d,-\e}$ $\ra G_{\d,\e}$. Then we can destabilize $\cC$ along $\r$.
\item[\rm(b)]  Suppose that as $\cC$ is traversed, a sub-arc $\r \subset \cC$ has the
following  sequential intersections with our four
graphs: $G_{\d,\e} \ra G_{-\d,-\e} \ra G_{\d,-\e} \ra G_{-\d,\e}$.
Then $\cC$ admits an exchange move along $\r$.
\end{itemize}
\end{lemma}

\pf
Referring \thinspace to \thinspace Figure \thinspace \ref{figure:standard tiling exchange move and
destabilization}(a), \thinspace we \thinspace see \thinspace
that \thinspace  a \thinspace sequential \thinspace intersection pattern
$G_{\d,\e} \ra G_{-\d,-\e} \ra G_{\d,-\e} \ra G_{-\d,-\e} \ra G_{\d,\e}$
implies the existence of a 
sub-disc in $\ta - \cC$ satisfying the assumptions in Lemma
\ref{lemma:destabilization}.  Conclusion (a) follows.

Referring again to Figure \ref{figure:standard tiling exchange move and
destabilization}(b), we observe 
that a sequential intersection pattern
$G_{\d,\e} \ra G_{-\d,-\e} \ra G_{\d,-\e} \ra G_{-\d,\e}$ implies the existence of a 
sub-disc in $\ta - \cC$ satisfying the assumptions in Lemma
\ref{lemma:ab-exchange move}.  Conclusion (b) follows. \endproof

\begin{figure}[ht!] 
\cl{\epsfysize=230pt \epsfbox{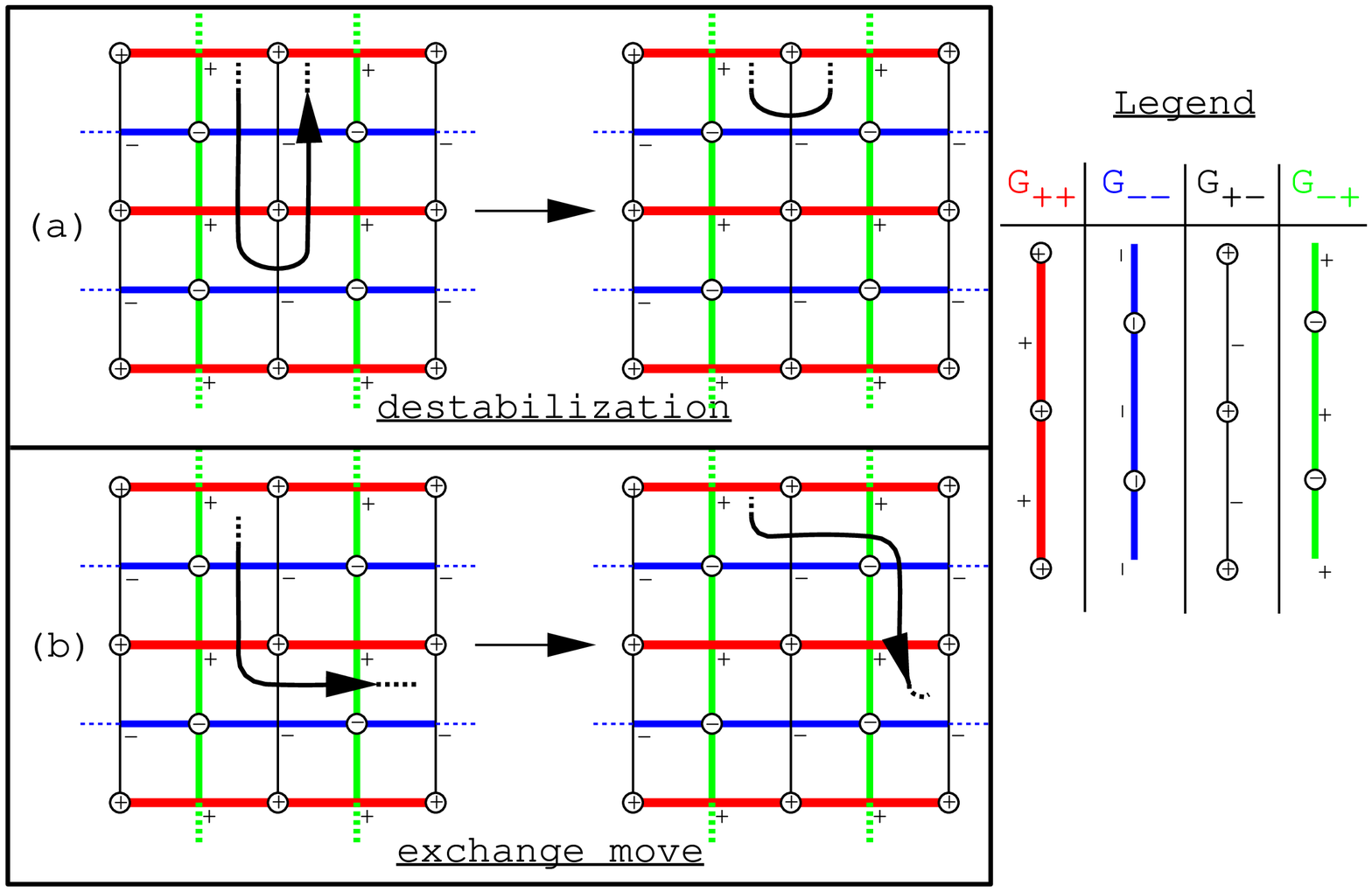}}
\caption{In (a) the intersection pattern is
$G_{-,-} \ra G_{+,+} \ra G_{-,+} \ra G_{+,+} \ra G_{-,-}$. In (b), the 
intersection pattern is $G_{-,-} \ra G_{+,+} \ra G_{+,-} \ra G_{-,+}$.}
\label{figure:standard tiling exchange move and destabilization}
\end{figure}

\begin{lemma}
\label{lemma:curve intersecting graphs on torus}
Let $\cC \subset \ta$ be a curve such that $\cC \in \{K,m\}$.
If $\cC$ intersects a component of $G_{\d,\e}$ incoherently then,
after a sequence of exchange moves, we can reduce the braid index of $\cC$ using
a destabilization.
\end{lemma}

\pf
Suppose that $\cC$ intersects $G_{+,+}$ incoherently.  Then there exists sub-arcs
$\g \subset G_{+,+}$ and $\r \subset \cC$ and a sub-disc $\Delta \subset \ta$
such that $\partial \Delta = \g \cup \r$.  Without loss of generality, we can assume
that $int(\Delta) \cap \cC = \emptyset$.  Now, notice that the standard tiling
of $\ta$ forces the existence of a sub-arc $\r^\prime \subset \r$ which satisfies the
sequential intersection pattern (a) or (b) of  Lemma
\ref{lemma:graph implies destabilization and exchange moves}.
Moreover, it can be assumed that the existence of the type (a) or (ab) vertex
from the proofs of Lemmas \ref{lemma:destabilization} and \ref{lemma:ab-exchange move}
(which the intersection sequences of Lemma
\ref{lemma:graph implies destabilization and exchange moves} invoke)
are contained in $\Delta$.

If the sequence in Lemma \ref{lemma:graph implies destabilization and exchange moves}(a)
occurs then we immediately have the conclusion of our lemma.
If the sequence in Lemma \ref{lemma:graph implies destabilization and exchange moves}(b)
occurs then, after performing the exchange move on $\r^\prime$, we will again
have sub-arcs $\g$ and $\r$.  But, they now bound a sub-disc $\Delta^\prime$
that has fewer vertices in its induced foliation than $\Delta$.  Iterating this
procedure we conclude that at some point $\r$ must have an intersection sequence
with the graphs that corresponds to the sequence in
Lemma \ref{lemma:graph implies destabilization and exchange moves}(a).
\endproof

We now focus on understanding how $K$ is contained in the foliation of $\ta$.
Specifically, we have the following application of Lemma \ref{lemma:curve intersecting graphs on torus}.

\begin{prop}
\label{proposition:K intersects graphs coherently}
Let $(K,\ta)$ be a knot--torus pair where $\ta$ has a standard tiling.
Then there exists a knot $K^\prime \subset \ta$ that is transverse to the foliation of
$\ta$ and is obtained from $K$ through a sequence of exchange moves and
destabilizations, with $\L(K^\prime,\ta) \leq \L(K,\ta)$ at every change,
such that $K^\prime$ coherently intersects components of the graphs
$G_{\pm,\pm}$ and $G_{\pm,\mp}$.  
\end{prop}

\pf
If $K$ does not intersect graph components coherently then, after repeated application
of Lemma \ref{lemma:curve intersecting graphs on torus}, we can assume that
$K$ has been replaced by $K^\prime$ satisfying the conclusion of the proposition.
If any of the $\bb$--arcs of $\ta$ are now inessential, we perform the necessary
isotopy of $\ta$ to remove them.  The new foliation of $\ta$ will have fewer
vertices and singularities.  We then repeat the applications of
Lemmas \ref{lemma:change of fibration} and \ref{lemma:eliminating valence two vertices},
emulating the argument at the begin of this section, so that $V(2)$ and $V(3)$
of equation (2) are again zero and the tiling of $\ta$ is checkerboard by the
singularity parity values.  All operations on $K$ and the tiling of $\ta$
are non-increasing on the complexity measure.
\endproof

For a scc $\cC \subset \ta$ transverse to the foliation of $\ta$, let
$ \cS_{\cC} \subset \ta$ be the closure of the union of all
the $\bb$--arcs that $\cC$ intersects in the foliation of $\ta$.   We call $\cS_{\cC}$ the
{\em $\bb$--support of $\cC$}. The definition of $\bb$--support implies a useful fact about
the ``width'' of $\cS_K$: If $\g\subset\cS_K$ is properly embedded (that is,
$(int(\g),
\partial\g)
\subset (int(\cS_K), \partial\cS_k)$, which is the union of $\bb$--arcs) then $\g$ is in
fact a single $\bb$--arc which is non-parallel to the boundary of $\cS_K$.  We will
refer to this fact as $(\star)$.

Next, we define two modifications of the $\bb$--support of $\cC$.

\begin{figure}[ht!] 
\cl{\epsfysize=160pt \epsfbox{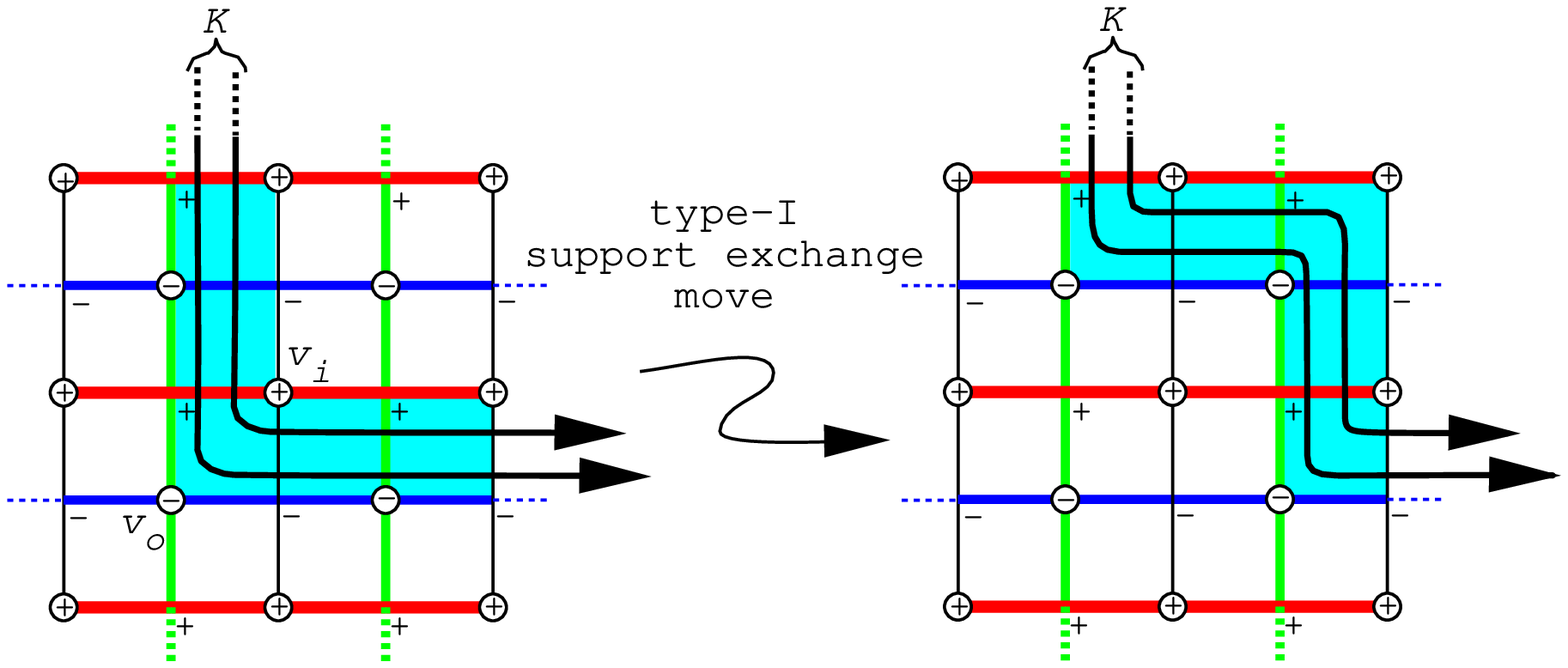}}
\caption{The alteration in $\cS_K$ from left to right corresponds to an
exchange move with type-I support.  The vertex labeled $v_o$ is an outside vertex and the
vertex labeled
$v_i$ is an inside vertex.  The shaded region is a portion of $\cS_K$.}
\label{figure:exchange move with type-I support}
\end{figure}

\noindent
{\bf Exchange move with type-I support}\qua Let $v_\e, v_{-\e} \subset \partial \cS_K$ be two vertices such that
if $b \subset \cS_K$ is a $\bb$--arc having $v_o$ as an endpoint then $v_i$ is also an endpoint of $b$.
More descriptively,  these conditions correspond to having $v_o$ as an {\em outside corner} and
$v_{i}$ is an {\em inside corner} of $\partial \cS_K$.  For
each point of $K \cap b$ we can find an arc neighborhood $b \subset \g \subset K$ such that $\g$ has the
sequential intersection pattern of
$ G_{-\e,-\d} \ra G_{\e,\d} \ra G_{\e,-\d} \ra G_{-\e,\d} $ as depicted in Figure
\ref{figure:standard tiling exchange move and destabilization}(b).  The alteration in
the $\bb$--support that results from performing the exchange move in Figure
\ref{figure:exchange move with type-I support} to each $\g$ arc intersecting $b$
is referred to as an {\em exchange move with type-I support}.

\begin{figure}[ht!] 
\cl{\epsfysize=180pt \epsfbox{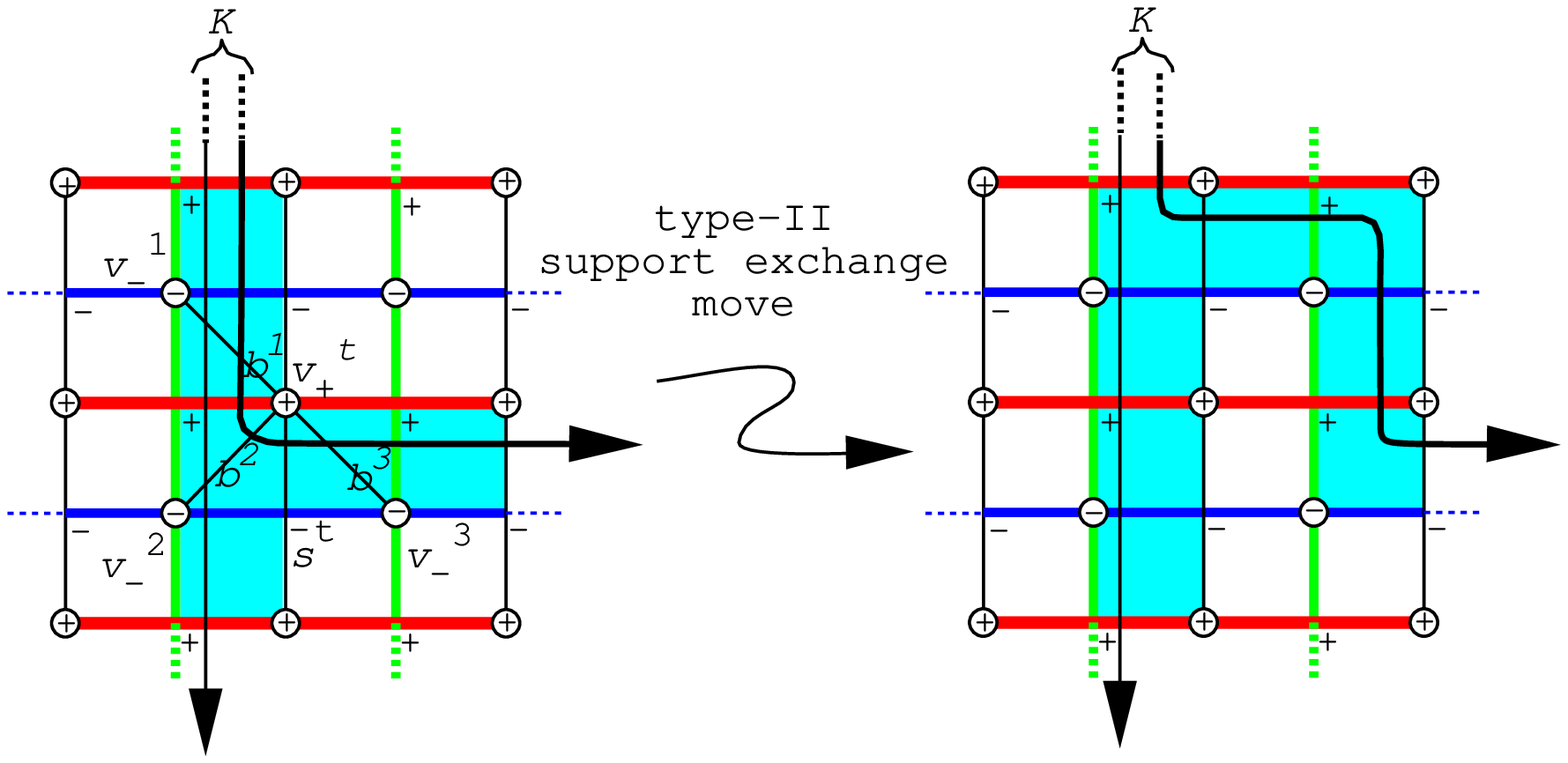}}
\caption{The alteration in $\cS_K$ from left to right corresponds 
to an exchange move with type-II support.  The vertex labeled $v_\e^t$ is a tee vertex
and the
 singularity labeled
$s^t$ is a tee singularity.  The shaded region indicates a portion of $\cS_K$.}
\label{figure:exchange move with type-II support}
\end{figure}

\noindent
{\bf Exchange move with type-II support}\qua  Let 
$v_\e^t,v_{-\e}^1,v_{-\e}^2,v_{-\e}^3 \subset \partial \cS_K$ be
four vertices such that if $b \subset \cS_K$ is a $\bb$--arc having 
$v_\e^t$ as an endpoint then its
other $b$ endpoint is contained in $\{v_{-\e}^1,v_{-\e}^2,v_{-\e}^3\}$.  
Let $b^1,b^2,b^3 \subset \partial \cS_K$
be three $\bb$--arcs such that $v_{-\e}^i$ is an endpoint of
$b^i$, $1 \leq i \leq 3$.  Assume the $i$ superscript
is in correspondence with the angular order of the $\bb$--arcs around 
the vertex $v_\e^t$.  Descriptively,
the vertex $v_\e^t$ is at a {\em tee vertex} of $\cS_K$.  The singular 
leaf that is common to vertices
$v_\e^t,v_{-\e}^1,v_{-\e}^2$ has its singular point on $\partial \cS_K$.  
This singular point, $s^t$, is a {\em tee
singularity} of $\cS_K$.
For any point of $K \cap b^2$ we can find an arc neighborhood
$b \subset \g \subset K$ such that either (path 1) $\g$ intersects in order 
$b^1, b^2, b^3$, or 
(path 2) $\g$ intersects in order
$b^1, b^2, b^\prime$ where $b^\prime$ is a $\bb$--arc having $v_{-\e}^2$ as 
an endpoint but not $v_{\e}$.
If only path 1 occurs then $v_{-\e}^2$ is a type of corner vertex of 
$\partial \cS_K$ and we have the situation
described in the exchange move with type-I support. So assume that both path 1 
and path 2 $\g$--arcs occur.
Since the sequential intersection pattern for path 1 $\g$--arcs corresponds 
to that of Figure
\ref{figure:standard tiling exchange move and destabilization}(b), we can 
perform a sequence of exchange moves
to $K$ until $v_\e$ no longer exists as a tee of $\cS_K$.  Figure 
\ref{figure:exchange move with type-II support}
shows the alternation to $\cS_K$ due to such a sequence of exchange moves.  
We refer to this alteration as an {\em exchange move with type-II support}.

A few final remarks about exchange moves with types I and II support will be useful.  
First, the four possible states of a fixed vertex on $\partial \cS_K$ are:
inside corner; outside corner; tee; and null, ie neither corner nor tee.
Second, suppose as we transverse a boundary component
$c \subset \partial\cS_K$ we encounter sequentially vertices $v^1,v^2,v^3$.  
Assume $v^2$ is an 
outside corner.  If
$v^1$(or $v^3$) has a null state then after an exchange move with type-I support at
$v^2$, we will have $v^1$ (or $v^3$) being an outside corner.
If $v^1$ (or $v^3$) is an inside corner then after an 
exchange move with type-I support at $v^2$ it has a null state. 
If $v^1$ (or $v^3$) is a tee
then after an exchange move with type-I support at $v^2$ it has a null state. Notice
if $K$ is coherent with respect to
$G_{\e,\d}$
neither $v^1$ nor $v^3$ may be an outside corner.
Now, assume $v_2$ is a tee vertex.
If $v^1$ and $v^3$ states are null then after an 
exchange move with type-II support at $v^2$, one vertex,
say $v^1$, will be a tee vertex and the other vertex, $v^3$, will be an 
inside corner vertex. If $v^1$ is an outside corner vertex then after 
an exchange move with type-II support at $v^2$, $v^1$
state will be null.  If there is a tee singularity being $v^1$ and $v^2$ then 
after
an exchange move with type-II support at $v^2$, $v^1$ will become a tee vertex.
Notice by the assumption that $K$ is coherent with respect to $G_{\e,\d}$ we 
cannot have
either $v^1$ or $v^3$ being an inside corner or tee vertex.

We now consider the possible subsurface types the $\bb$--support might assume.

\begin{lemma}
\label{lemma:annulus b-support must be homotopic to G_{+,+} components}
Let $(K,\ta)$ be a knot--torus pair where $\ta$ has a standard tiling.
Assume that $K$ intersects all components of the graphs
$G_{\pm,\pm}$ and $G_{\pm,\mp}$ coherently and that the $\bb$--support of $K$ is an
annulus in $\ta$.  If $\partial \cS_k$ contains a corner then there exists a
$(K^\prime,\ta^\prime)$ that is obtained from $(K,\ta)$ by a sequence of exchange moves
such that $\L(K^\prime,\ta^\prime) < \L(K,\ta)$.
\end{lemma}

\pf
First, notice that if $v \subset \partial \cS_K$ is an outside or inside corner 
and the unique $\bb$--arc
in $\cS_K$ that has $v$ as its endpoint is an outermost essential arc, then, 
after an exchange move with type-I support at $v$,
this $\bb$--arc will become inessential.

Second, suppose that as we transverse a boundary component $c \subset \partial\cS_K$ 
we encounter
sequentially vertices $v^1,v^2,\cdots,v^l$, where $v_l$ is the only corner vertex in our list.
Furthermore, suppose that $v_1$ is adjacent to an outermost essential $\bb$--arc, 
$b \subset \cS_K$.
We can then ``walk'' the corner at $v_l$ up to $v_1$: perform 
an exchange move with type-I support  which
will result in $v_{l-1}$ being a corner; and iterate this process along $c$ until 
$v_1$ is a corner.
Finally, we perform an exchange move with type-I support at $v_1$ to make $b$ inessential.

Removing the inessential $\bb$--arc will produce a new knot--torus pair having 
decreased complexity.
\endproof

\begin{lemma}
\label{lemma:scc in boundary of support bounding a disc must have a tee or corner}
Let $(K,\ta)$ be a knot--torus pair where $\ta$ has a standard tiling.
Assume that $K$ intersects all components of the graphs
$G_{\pm,\pm}$ and $G_{\pm,\mp}$ coherently and that $c \subset \partial \cS_K$ is a scc
which bounds a disc
$\Delta_c \subset \ta - int(\cS_K)$.  If $|\Delta_c|_v$ is the number of vertices in the interior
of $\Delta_c$ then there exists an exchange move with special support that will reduce
$|\Delta_c|_v$.
\end{lemma}

\pf
Since $c$ is the union of edge-paths in the graphs $G_{\e,\d}$ and since the 
tiling of $\ta$ is standard (specifically,
all vertices are valence four), we know that $c$ must contain either an 
inside corner or a tee vertex.
If $c$ contains an inside corner then after an exchange move with type-I support, 
the number of vertices contained
in the interior of $\Delta_c$ will decrease.  If $c$ contains a tee vertex 
then after an 
exchange move with type-II support, again, the number of vertices contained in the 
interior of $\Delta_c$ will decrease.
\endproof

\begin{prop}
\label{proposition:the b-support of K}
Let $(K,\ta)$ be a knot--torus pair where $\ta$ has a standard tiling
and assume that $K$ coherently intersects all components of the graphs
$G_{\pm,\pm}$ and $G_{\pm,\mp}$.
Then there exists a knot $K^\prime \subset \ta$ that is transverse to the foliation of
$\ta$ and is obtained from $K$ through a sequence of exchange moves, with
$\L(K^\prime,\ta^\prime) = \L(K,\ta)$,
such that either:
\begin{itemize}
\item[\rm(i)]  $\cS_{K^\prime}$ is an annulus with $\partial\cS_{K^\prime} = c_1 \cup c_2$
where
$c_1$ is a component of $G_{\e,\d}$ and $c_2$ is a component of $G_{-\e,-\d}$;
\item[\rm(ii)] $\cS_{K^\prime}$ is a torus-minus-a-disc with $\partial\cS_{K^\prime} = \a_1 \cup \b_1 \cup \a_2 \cup \b_2$
where $\a_1$ and $\a_2$ are arcs in $G_{\e,\d}$ and $\b_1$ and $\b_2$ are arcs in $G_{\e,-\d}$.  In particular,
$\partial \cS_K$ does not contain any corner vertices.
\end{itemize}
\end{prop}

\pf
Since $\cS_{K} \subset \ta$ and $\ta - K$ is an annulus, we can distinguish the 
cases for $\cS_K$ topologically as being
either an annulus, an annulus minus discs, or a torus minus discs. 
\begin{itemize}
\item[(i-a)]
If $\cS_K$ is an annulus then,
by Lemma \ref{lemma:annulus b-support must be homotopic to G_{+,+} components} 
we have the annulus in (i).
\item[(i-b)] If $\cS_K$ is an annulus minus discs let $c \subset \partial \cS_K$
be a component that bounds a disc $\Delta_c \subset \ta - int(\cS_K)$.
By Lemma \ref{lemma:scc in boundary of support bounding a disc must have a tee or corner}
we can reduce the number of vertices contained in $int(\Delta_c)$ iteratively until
we eliminate a component
of $\partial \cS_K$.  Thus, an annulus minus discs can be reduced to an annulus.  
Appealing to
Lemma \ref{lemma:annulus b-support must be homotopic to G_{+,+} components} again, 
we reduce to the annulus in (i).
\item[(ii)] If $\cS_K$ is a torus minus discs, we can apply the
argument above to reduce the number of components of $\partial \cS_K$ to one component. 
If this single boundary component contains a corner then, as in the proof of Lemma
\ref{lemma:annulus b-support must be homotopic to G_{+,+} components}, we can ``walk''
that corner past any outermost essential $\bb$--arc in $\cS_K$.  So $\cS_K$ will be a
torus minus a disc as described in (ii). \qed
\end{itemize}

\subsection{The knot, torus, meridian triple}
\label{subsection:knot, torus, meridian triple}
We now consider the triple $(K,\ta,m)$ where $m \subset \ta$ is a 
meridian curve that intersects $K$ coherently
$p_h$--times.  We expand our measure of complexity to be the $3$--tuple 
$\L(K,\ta,m) = (n_1 , n_2 , n_3)$
where $n_3$ is the braid index of $m$ and, as before,
$n_1$ is the braid index of $K$ and 
$n_2$ is the number of vertices in the foliation of $\ta$.

\begin{prop}
\label{proposition:m intersects graphs coherently}Let $(K,\ta,m)$ be a knot--torus--meridian triple where the pair $(K,\ta)$ 
satisfies the conclusion of Proposition \ref{proposition:K intersects graphs coherently}.
Then there exists a meridian curve $m^\prime \subset \ta$ that is transverse to the foliation of
$\ta$ and is obtained from $m$ through a sequence of exchange moves and destabilization
such that $m^\prime$ intersects all components of the graphs
$G_{\pm,\pm}$ and $G_{\pm,\mp}$ coherently .  Moreover, $\L(K,\ta,m^\prime) \leq
\L(K,\ta,m)$.
\end{prop}

\pf
We model our proof on the proof of Proposition \ref{proposition:K intersects graphs
coherently}.  If $m$ does not intersect graph components coherently then, after repeated
application of Lemma \ref{lemma:curve intersecting graphs on torus}, we can assume that
$m$ has been replaced by $m^\prime$ satisfying the conclusion of the proposition.
Notice that if $m$ is destabilized using the isotopy in Figure
\ref{figure:standard tiling exchange move and destabilization}, this isotopy will
occur away from $K$. Moreover,
if $m$ is isotopied through a Figure \ref{figure:standard tiling exchange move and
destabilization} exchange move then $m$ will still intersect $K$ $p_h$--times.
All operations on $m$
are non-increasing on the complexity measure.
\endproof

Let $v_+$ and $v_-$ be two vertices in the foliation of $\ta$ which are endpoints of a common $\bb$--arc.
Consider the rectangular region in $\ta$ that is the closure of the union of all the $\bb$--arcs that
have $v_+$ and $v_-$ as their endpoints.  We shall call such a region a 
{\em $\bb$--rectangle}.

\begin{figure}[ht!] 
\cl{\epsfysize=130pt \epsfbox{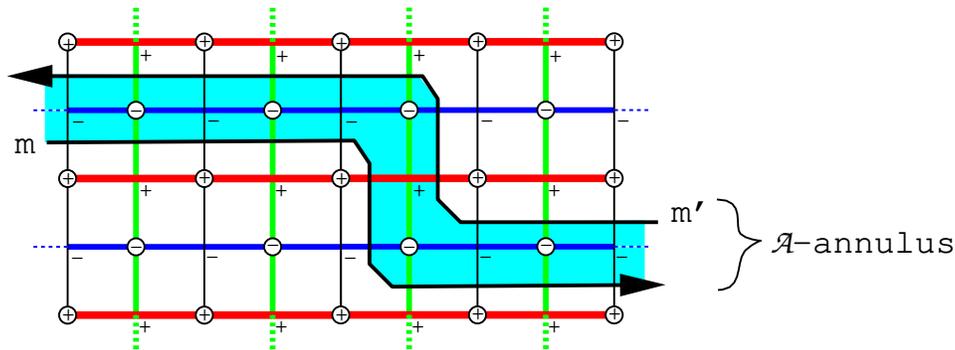}}
\caption{The shaded region between $m$ and $m^\prime$ is $\cA$.  Notice that the
induced foliation of $\cA$ has only $\ba\ba$--singularities.}
\label{figure:m is isotopic to mprime}
\end{figure}

We introduce the notation of a {\em parallel push-off} of a meridian curve $m$.
Specifically, we can choose a scc $m^\prime \subset \ta$ such that:
there exists an annulus $\cA \subset \ta$ with $\partial \cA = m \cup m^\prime$; the induced foliation
on $\cA$ has only $\ba\ba$--singularities; and each vertex in the foliation of $\cA$ is adjacent to two
singular leaves.  See Figure \ref{figure:m is isotopic to mprime}.  The oriented
arc $m$ is isotopic on $\ta$ to the oppositely oriented arc $m^\prime$.
We use the notation $\cS_m$ for the $\bb$--support of the meridian curve $m$.

\begin{prop}
\label{proposition:support of K intersect the support of m}
Let $(K,\ta,m)$ be a knot--torus--meridian triple satisfying
the conclusions of Propositions \ref{proposition:K intersects graphs coherently},
\ref{proposition:the b-support of K} and
\ref{proposition:m intersects graphs coherently}.
Then we can assume that $\cS_K \cap \cS_m$ is a union of disjoint $\bb$--rectangles.
\end{prop}

\pf
By Lemma \ref{lemma:annulus b-support must be homotopic to G_{+,+} components},
we know that $\cS_K$ is either an annulus or a torus minus a disc.  Moreover, we
also know by the same lemma that $\partial \cS_K$ has no corners.

Using the same argument that was used in the proof of Lemma
\ref{lemma:annulus b-support must be homotopic to G_{+,+} components}, we can assume that each
boundary component of $\cS_m$ has at most one outside (inside) corner.

If $\cS_K$ is an annulus of the type-(ii) in Proposition \ref{proposition:the b-support of K}
and $\partial\cS_m$ has no corners
then we immediately have the conclusion of the proposition.  If $\partial\cS_m$ has any corners then
we can ``walk'' these corners
(as in the proof of Lemma \ref{lemma:annulus b-support must be homotopic to G_{+,+} components})
so that the corners of $\partial\cS_m$ are
away from $\cS_K \cap \cS_m$.  After this repositioning, $\cS_K \cap \cS_m$ will be disjoint $\bb$--rectangles.

It is easily seen that if $\cS_K \cup \cS_m$ is not a disjoint union of
$\bb$--rectangles then $\cS_K \cup \cS_{m^\prime}$ will be a disjoint union of $\bb$--rectangles
where $m^\prime$ is a parallel push-off of $m$.
\endproof

\subsection{The knot, torus, meridian and meridian-disc}
\label{subsection:knot, torus, meridian, and meridian-disc}
We now consider the quadruple
$(K,\ta,m,\Delta_m)$ where $\Delta_m$ is a meridian disc which $m$ bounds inside the
solid torus $\va$ which $\ta$ bounds. Recall that the fibration $\fibr$ induces a
foliation on
$\Delta_m$. Again, we expand our measure of complexity to be the $4$--tuple 
$\L(K,\ta,m,\Delta_m)
= (n_1 , n_2 , n_3, n_4)$ where $n_4 = |\axis \cap \Delta_m|$ and $n_1$, $n_2$ and
$n_3$ are as before.  If $n_4=0$ then the foliation of $\ta$ contains $\bc$--circles and we
can appeal to the arguments in Section \ref{mixed foliations}.

Before we proceed further it will be useful to review the induced foliation on
$\Delta_m$.  Again, a comprehensive reference is \cite{[BF]}.  When $n_4 = 1$,
$\Delta_m$ is {\em radially} foliated by
$\ba$--arcs adjacent to the unique vertex.  When $n_4 >1$, the singular foliation of $\Delta$ contains
the vertices $\Delta_m \cap \axis$ and
saddle singularities of three possible types---$\ba\ba$--, $\ba\bb$-- and $\bb\bb$--singularities.
Thus, $\Delta_m$ is tiled by $\ba\ba$, $\ba\bb$ and $\bb\bb$ tiles.
Since $m$ induces on $\Delta_m$ an orientation, we can
assign a parity to each vertex and singularity using the same assignment scheme employed in the
foliation of $\ta$.  Having a parity assignment on vertices and singularities allows us to
define graphs $G_{\pm,\pm}$ and $G_{\pm, \mp}$ as before with the proviso that we treat
$m$ as if it were a negative vertex.  (This implies that $G_{-,-}$ and $G_{-,+}$ will contain
edges having endpoints on $m$.)

\begin{figure}[ht!] 
\cl{\epsfxsize=300pt \epsfbox{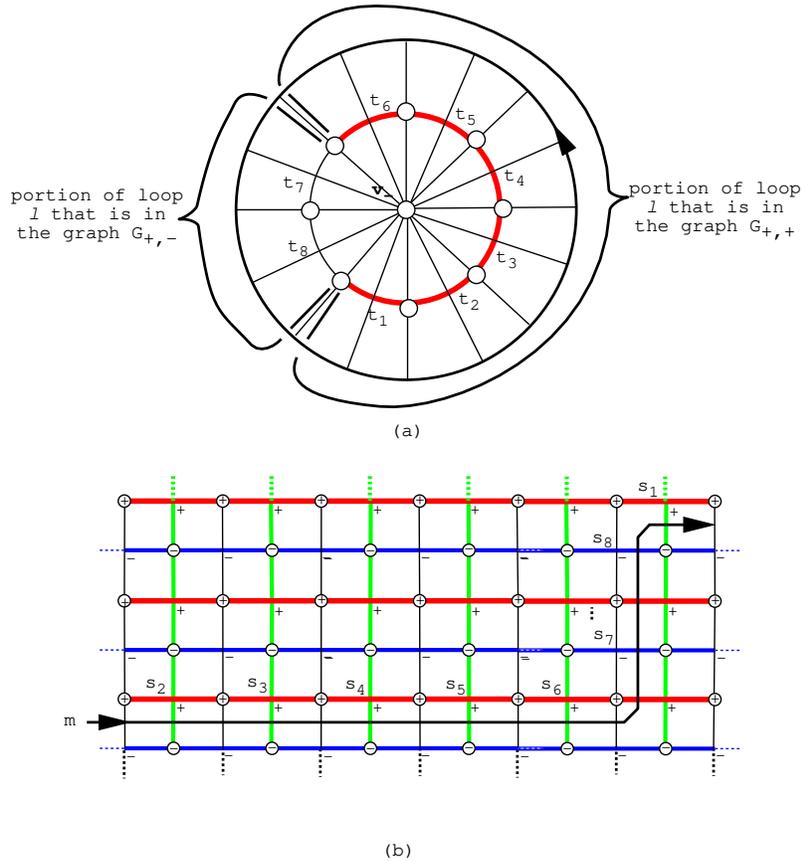}}
\caption{In (a) we have the foliation of $\Delta_m$ 
illustrated in the special case
where $l$ contains two negative and six positive singularities, and there is a single negative
vertex.  In (b) we have the
corresponding path $m$ in a standard tiling of $\ta$.
The arrowheads are
identified to form a scc.  The labels $\{t_1, \cdots, \t_8 \}$ of the 
singularities in (a)
show the sequential correspondence to the singularities in (b) with the labels
$\{s_1, \cdots , s_8 \}$.}
\label{figure:initial condition on Deltam}
\end{figure}

Now consider how the foliations of $\Delta_m$ and $\Delta_{m^\prime}$ are
related to one another, when $m^\prime$ is a parallel push-off of $m$.  
Let $\cA \subset \ta$ be the 
annulus in Figure \ref{figure:m is isotopic to mprime} that has $m$ and 
$m^\prime$ as its
boundary curves.  Notice that the core circle in $\cA$ that is contained 
in a union of all of
the singular leaves in the induced foliation of $\cA$ is a union of sub-arcs
 of the graphs
$G_{\e,\d} \subset \ta$.  Similar in flavor to the arguments in Proposition 
\ref{proposition:the b-support of K},
we can produce a sequence of exchange moves by applying Lemma
\ref{lemma:graph implies destabilization and exchange moves} so that either
 $\partial\cS_m$ contains
no corners, or each component of $\partial\cS_m$ contains exactly one inside 
corner and one outside
corner.  If $\partial\cS_m$ contains no corners then the foliation of $\cA$ will have singularities
all of the same parity.  If each component of $\partial\cS_m$ has exactly one inside and outside corner
then the core circle in $\cA$ is the union of an arc in $G_{\e,+}$ and an arc in $G_{\e,-}$.

Now the foliation of $\Delta_{m^\prime}$ is simply the foliation of $\Delta_m
\cup_m \cA$. Specifically, for each singularity in $s \in \cA \subset \ta$ there is
a singularity
$t \in G_{+,\pm} \subset \Delta_{m^\prime}$ (coming from the ``$\cA$ portion'' of $\Delta_m$)
that is in a regular neighborhood of $s$ on the negative side of the
oriented $\ta$.  It is easily checked by appealing to the $H_t$--sequence of $\ta \cup \Delta_m \subset \fibr$ that
if $s$ is a positive (negative)
singularity then $t$ is the singularity that occurs immediately before (after)
$s$. (See Figure \ref{figure:initial condition on Deltam}.)
If $\partial\cS_m$ contains no corners then the $G_{\e,\d}$ of
$\Delta_{m^\prime}$ contains
a loop---the core circle coming from $\cA$.  But by Lemma 3.8 of \cite{[BF]}
such a loop cannot exist.  So we are left having $\Delta_{m^\prime}$
 containing a loop (which is
the core circle of $\cA \subset \Delta_{m^\prime}$) that is the union 
of an arc in the graph
$G_{+,+} \subset \Delta_{m^\prime}$ and an arc in the graph 
$G_{-,-} \subset \Delta_{m^\prime}$.
This loop contains all of the positive vertices of $\Delta_{m^\prime}$.
As noted above, a singularity $s \in \partial\cS_m$ will in $\fibr$
sequentially correspond to a singularity $t \in l$.
An application of the argument in Lemma 3.8 (ii) allows us to conclude 
this discussion of the foliation
of $\Delta_{m^\prime}$ with the following summarizing result.

\begin{lemma}
\label{lemma:the intial foliation of Delta_m}
Let $(K,\ta,m,\Delta_m)$ be a quadruple where the triple $(K,\ta,m)$ 
satisfies the conclusions of Propositions \ref{proposition:K intersects graphs coherently},
\ref{proposition:the b-support of K}, \ref{proposition:m intersects graphs coherently} and
\ref{proposition:support of K intersect the support of m}.
We can assume that the initial foliation of $\Delta_m$ satisfies the following conditions.
\begin{itemize}
\item[\rm(a)] There are only $\ba\bb$--singularities.
\item[\rm(b)] There exists a closed loop $l$ that is the union of an arc in
$G_{+,+} \subset \Delta_m$ and an arc in $G_{+,-} \subset \Delta_m$.
\item[\rm(c)] When $n_4=3$, there is a single negative vertex $v_-$.
When $n_4 > 3$, there is a negative vertex $v_-$ adjacent to two singularities of common parity that are
consecutive on $l$.
\end{itemize}
Moreover, the
singularities contained in each component of $\partial\cS_m$ are in one-to-one correspondence
with, and in $\fibr$ sequentially correspond to, the singularities contained in $l$.
(See Figure \ref{figure:initial condition on Deltam}.)
\end{lemma}

We now use this initial tiling of $\Delta_m$ to reduce the complexity of the quadruple
$(K,\ta,m,\Delta_m)$.  Notice that the number of singularities in the loop
$l \subset \Delta_m$ is naturally two or greater.  The next lemma shows us that
when $l$ contains exactly two singularities the foliation of a tiled $\ta$ must have
an inessential $\bb$--arc.

\begin{lemma}
\label{lemma:Delta_m has n4=3}
Let $(K,\ta,m,\Delta_m)$ be a quadruple, where the triple $(K,\ta,m)$ 
satisfies the conclusions of Propositions \ref{proposition:K intersects graphs coherently},
\ref{proposition:the b-support of K}, \ref{proposition:m intersects graphs coherently},
\ref{proposition:support of K intersect the support of m} and
\ref{lemma:the intial foliation of Delta_m}.
If $n_4 =3$ then the foliation of $\ta$ has an inessential $\bb$--arc.
\end{lemma}

\pf
If $n_4 = 3$ then $l$ must have contained one positive and one negative singularity.
Thus, referring to the terminology of \cite{[BM3]}, $\ta$ must be a type {\bf k} torus.
A type {\bf k} torus embedding has the property that a meridian curve on $\ta$ can be
represented by the union of two $\bb$--arcs, each being outermost in a disc fiber of
$\fibr$.  Since
$K$ can only intersect one of these two $\bb$--arcs, $\ta$ must contain an inessential $\bb$--arc.
(Refer back to Remark \ref{remark:why the proof is difficult} and Figure \ref{figure:example of a tiled torus}
for an example of a type {\bf k} torus.)
\endproof

We next deal with the case where $l$ contains more than two singularities.

\begin{lemma}
\label{lemma:Delta_m has n4>3}
Let $(K,\ta,m,\Delta_m)$ be a quadruple where the triple $(K,\ta,m)$ 
satisfies the conclusions of Propositions \ref{proposition:K intersects graphs coherently},
\ref{proposition:the b-support of K}, \ref{proposition:m intersects graphs coherently},
\ref{proposition:support of K intersect the support of m} and
\ref{lemma:the intial foliation of Delta_m}.
If $n_4 > 3$ then there exists an isotopy of $\ta$ such that $n_2$ is reduced.
\end{lemma}

\begin{figure}[ht!] 
\cl{\epsfxsize=350pt \epsfbox{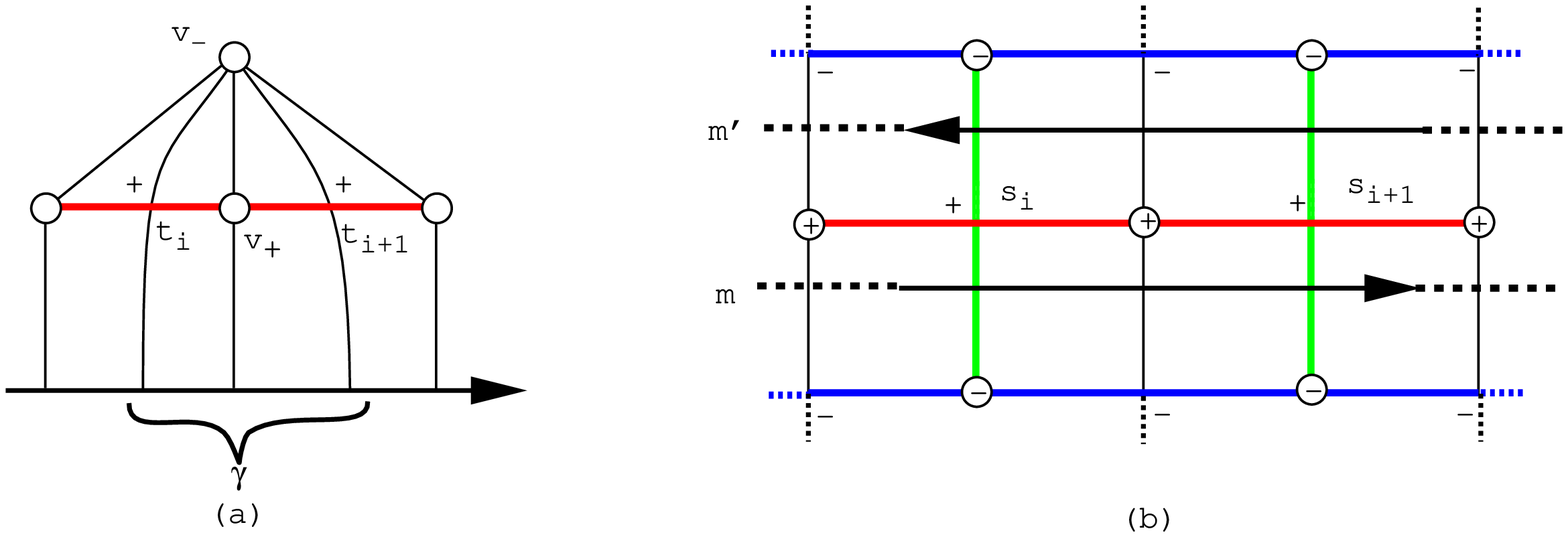}}
\nocolon \caption{}
\label{figure:isotopy of ta across valence 1 vertex part 1}
\end{figure}

\pf
If $n_4 >3$ then we employ the negative vertex $v_-$ in condition (c) of Lemma
\ref{lemma:the intial foliation of Delta_m}.  Figure
\ref{figure:isotopy of ta across valence 1 vertex part 1}(a) illustrates this situation
when the parity of these two singularities is positive.  In Figure
\ref{figure:isotopy of ta across valence 1 vertex part 1}(a) these singularities 
are labeled
$t_{i}$ and $t_{i+1}$. The common positive vertex that is adjacent to both
of their associated singular leaves is $v_+$.  A regular neighborhood of the portion of
$m$ that contains endpoints of $\ba$--arcs adjacent to $v_+$ is labeled $\g$.  If we look
at the image of $\g$ in the foliation of $\ta$ we see that $\g$ intersects  a
positive singular leaf containing the singular point $s_i$ associated to
$t_i$; then a negative singular leaf; then another positive singular leaf containing
the singular point $s_{i+1}$ associated to $t_{i+1}$.  In the $H_t$--sequence
of $\Delta_m \cup \ta$ we know by the discussion in the proof of Lemma
\ref{lemma:the intial foliation of Delta_m} that $t_i$ occurs immediately before $s_i$
and $t_{i+1}$ occurs immediately before $s_{i+1}$.

By Proposition \ref{proposition:support of K intersect the support of
m}, there are now three possibilities for the arc $\gamma$: (1) $\g \cap K =
\emptyset$; (2) $\g$ intersects the singular leaf belonging to $s_i$ first, then the
negative singular leaf and then $K$;  or (3) $\g$ intersects
the singular leaf belonging to $s_i$  first, then intersects $K$ some number of times,
and then intersects the negative singular leaf.

If (3) occurs we replace $m$ with $m^\prime$, altering $\Delta_{m^\prime}$
so that its foliation again satisfies the initial condition foliation of
Lemma \ref{lemma:the intial foliation of Delta_m}.  This produces a new $\g$
which corresponds to possibility (2).  Now notice in the situation of (2)
we can isotopy $K$ in $\ta$ so
that $K$ intersects $m$ between the time $m$ crosses the $t_{i+1}$ singular leaf and the
time it crosses the $s_{i+1}$ singular leaf.  Thus, we can assume that $\g \cap K =
\emptyset$.
\begin{figure}[ht!] 
\cl{\epsfxsize=350pt \epsfbox{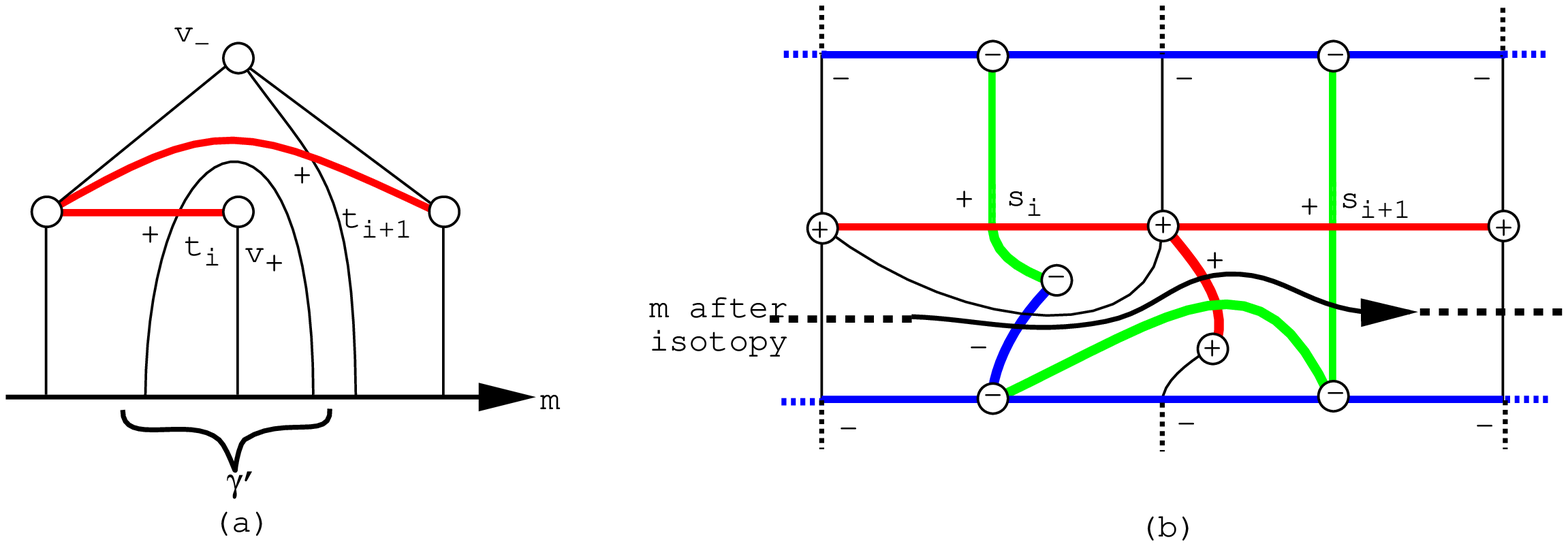}}
\nocolon \caption{}
\label{figure:isotopy of ta across valence 1 vertex part 2}
\end{figure}

Since $t_i$ and $t_{i+1}$ are of the same parity and adjacent to the unique
negative vertex $v_-$ in $\Delta_m$, we can perform a change of foliation.  (This
is an application to the foliation of $\Delta_m$ of the alteration in Figure
\ref{figure:change of fibration}.) Figure \ref{figure:isotopy of ta across valence 1
vertex part 2}(a) show how this change of foliation results in a type (a) vertex, $v_+
\subset \Delta_m$, which locally has the same foliation as that in Figure
\ref{figure:destabilization}.  The boundary of this valence $1$ vertex contains
$\g^\prime \subset m$.  (See Figure
\ref{figure:isotopy of ta across valence 1 vertex part 2}(a).)  We now destabilize $\Delta_m$ through
this valence $1$ vertex and, since $\g^\prime \cap K = \emptyset$, we can drag $\ta$ along through
this destabilization of $m$ without altering the embedding of $K$.  Figure
\ref{figure:isotopy of ta across valence 1 vertex part 2}(b) illustrates how the foliation of
$\ta$ is altered by this destabilization of $m$.  Notice that two new vertices are introduced into
the foliation of $\ta$, one positive and one negative, and that both of these new
vertices are valence $2$.   We can now appeal to
Figure \ref{figure:valence two vertex} to eliminate four vertices of $\ta$.
This reduces $n_2$ and, thus, the complexity of our quadruple.

Finally, \ it \ should \ be \ noted \ that \ we \ have made a choice of the parity of
$\{t_i,t_{i+1},s_i,s_{i+1}\}$ for reasons having to do with the clarity of the
expository. This does not reduce the generality of the argument. \endproof

We have proved Proposition \ref{prop:final goal} and, thus, Theorem \ref{calculus on iterated torus knots}.

\end{document}